\documentclass[10pt,leqno]{article}

\usepackage{amsmath,amssymb,amsthm,mathrsfs,dsfont}

\usepackage[margin=3cm]{geometry}

\usepackage{titlesec,hyperref}

\usepackage{color}

\usepackage{fancyhdr}
\titleformat{\subsection}{\it}{\thesubsection.\enspace}{1pt}{}

\newtheorem{theo}{Theorem}[section]
\newtheorem{lemm}[theo]{Lemma}
\newtheorem{defi}[theo]{Definition}
\newtheorem{coro}[theo]{Corollary}
\newtheorem{prop}[theo]{Proposition}
\newtheorem{rema}[theo]{Remark}
\numberwithin{equation}{section}

\allowdisplaybreaks

\newcommand\ep{{\varepsilon}}
\newcommand\lm{{\lesssim}}

\begin{document}
\title{Global existence and well-posedness for the Doi-Edwards polymer model
\hspace{-4mm}
}

\author{Wei $\mbox{Luo}^1$\footnote{E-mail:  luowei23@mail2.sysu.edu.cn} \quad and\quad
 Zhaoyang $\mbox{Yin}^{1,2}$\footnote{E-mail: mcsyzy@mail.sysu.edu.cn}\\
 $^1\mbox{Department}$ of Mathematics,
Sun Yat-sen University,\\ Guangzhou, 510275, China\\
$^2\mbox{Faculty}$ of Information Technology,\\ Macau University of Science and Technology, Macau}
\date{}
\maketitle
\hrule

\begin{abstract}
In this paper we mainly investigate the Cauchy problem of the Doi-Edwards polymer model with dimension $d\geq2$. The model was derived in the late 1970s to describe the dynamics of polymers in melts. The system contain a Navier-Stokes equation with an additional stress tensor which depend on the deformation gradient tensor and the memory function. The deformation gradient tensor satisfies a transport equation and the memory function satisfies a degenerate parabolic equation. We first proved the local well-posedness for the Doi-Edwards polymer model in Besov spaces by using the Littlewood-Paley theory. Moreover, if the initial velocity and the initial memory is small enough, we obtain a global existence result. \\

\vspace*{5pt}

\vspace*{5pt}
\end{abstract}

\vspace*{10pt}


\tableofcontents

\section{Introduction}
   In this paper we consider the Doi-Edwards polymer model \cite{Doi1980,Doi1988}:
   \begin{align}\label{1.1}
\left\{
\begin{array}{ll}
Re [\partial_{t}v+(v\cdot\nabla)v]-(1-\omega)\Delta{v}+\nabla{P}=div~\sigma, ~~~~~~~div ~v=0,\\[1ex]
\sigma(t,x)=\omega\displaystyle\int^{\frac{1}{2}}_{-\frac{1}{2}}S(t,x,s)ds,~~S(t,x,s)=-\displaystyle\int^{+\infty}_0\partial_{\tau}K(t,\tau,x,s)\mathscr{S}(G(t,\tau,x))d\tau,\\[1ex]
\partial_{t}G+(v\cdot\nabla)G+\frac{1}{We}\partial_\tau G=G\cdot\nabla v,  \\[1ex]
\partial_tK+(v\cdot\nabla)K+\frac{1}{We}\partial_\tau K+(\nabla v:\int^s_0S)\partial_sK-\frac{1}{We}\partial^2_sK=0, \\[1ex]
v(t,x)|_{t=0}=v_{0},\quad G(t,\tau,x)|_{t=0}=G_{0}(\tau,x) ,\quad K(t,\tau,x,s)|_{t=0}=K_0(\tau,x,s), \\[1ex]
K(t,\tau,x,s)|_{\tau=0}=1,\quad K(t,\tau,x,s)|_{s=\frac{1}{2}}=K(t,\tau,x,s)|_{s=-\frac{1}{2}}=0, \quad  G(t,\tau,x)|_{\tau=0}=I.\\[1ex]
\end{array}
\right.
\end{align}

In fact, the system (\ref{1.1}) is not the original model but a dimensionless model \cite{Chupin2017}. In $(1.1)$ $v$ stands for the velocity of the polymeric liquid which depends on the current time $t$ and the spatial variable $x\in\mathbb{R}^{d}$ with $d\geq2$. The deformation gradient tensor $G$ depends not only on the current time $t$ and the spatial variable $x$ but also on the history time $\tau$. The integral kernel $K$ depends on $t,\tau,x$ and the arc-length quantity $s\in[-\frac{1}{2},\frac{1}{2}]$. $\sigma$ is the stress tensor which can be deduced from a tensor $S$. The tensor $S$ can be written into a integral form. $Re$ is the Reynolds number and $We$ is the Weissenberg number. There are two viscosity in a dilute polymer flow: the solvent viscosities $\eta_s$ and the elastic viscosity $\eta_e$. The viscosities ratio $w=\frac{\eta_e}{\eta}\in(0,1)$ where $\eta=\eta_e+\eta_s$ is the total viscosity. Finally, the model is closed with the expression of the function $\mathscr{S}$:
\begin{align}
\mathscr{S}(G)=\frac{1}{<|G\cdot R|>_0}\langle\frac{(G\cdot R)\otimes(G\cdot R)}{|G\cdot R|}\rangle_0-\frac{I}{d},
\end{align}
where $<f>_0=\int_{\mathbb{S}^{d-1}}fdR$ and $R$ is a unit vector on $\mathbb{S}^{d-1}$.

   This model describes the dynamics of flexible polymers in melts. The system is of great interest in many branches of physics, chemistry, and biology, see \cite{Doi1988}. Roughly speaking, the system couples the incompressible Navier-Stokes equations with a transport system describing the deformation of polymers and a transport-diffusion equation describing the evolution of the polymer memory.

   A famous model to describe a viscoelastic flow is the Oldroyd-B model as follows:
   \begin{equation}
 \left\{\begin{array}{ll}
  \partial_tv+v\cdot\nabla v-\mu\Delta v+\nabla P=div~(G^TG),\\
  \partial_tG+v\cdot\nabla G= G\cdot\nabla v,\\
  div~v=div~G=0, \\
  v(t,x)|_{t=0}=v_0,~~G(t,x)|_{t=0}=G_0.
  \end{array}\right.
  \end{equation}
This model has been extensively studied. The equation $div G=0$ can be deduce from the condition on initial data such that $div G|_{t=0}=0$, thus the model is closed. Let us review some mathematical results for the Oldroryd-B model. C. Guillop\'e and J.C. Saut \cite{Guillope1990-NA,Guillope1990} proved the existence of local strong solutions and the global existence of one dimensional shear flows. In \cite{Fernandez-Cara}, E. Fern\'andez-Cara, F. Guill\'en and R. Ortega studied the local well-posedness in Sobolev spaces. J. Chemin and N. Masmoudi \cite{Chemin2001} proved the local well-posedness in critical Besov spaces and give a low bound for the lifespan. P. L. Lions and N. Masmoudi \cite{Lions-Masmoudi} obtained a global existence result of weak solutions. In \cite{Lin-Liu-Zhang2005}, F. Lin, C. Liu and P. Zhang proved that if the initial data is a small perturbation around equilibrium, then the strong solution is global in time. The similar results were obtained in several papers by virtue of different methods, see Z. Lei and Y. Zhou \cite{Lei-Zhou2005}, Z. Lei, C. Liu and Y. Zhou \cite{Lei2008}, T. Zhang and D. Fang \cite{Zhang-Fang2012}. For the Oldoryd-B model, the global existence of strong solutions in two dimension without small conditions is still an open problem.

The other famous model to describe the dynamic of polymer flow is the finite extensible nonlinear elastic (FENE) dumbbell model \cite{Bird1977,Doi1988}:
   \begin{align}
\left\{
\begin{array}{ll}
\partial_{t}v+(v\cdot\nabla)v-\mu\Delta{v}+\nabla{P}=div~\sigma, ~~~~~~~div~v=0,\\[1ex]
\partial_{t}\psi+(u\cdot\nabla)\psi=div_{R}[-\nabla{u}\cdot{R}\psi+\beta\nabla_{R}\psi+\nabla_{R}\mathcal{U}\psi],  \\[1ex]
\sigma_{ij}=\int_{B}(R_{i}\nabla_{j}\mathcal{U})\psi dR. \\[1ex]
\end{array}
\right.
\end{align}

This system is a micro-macro model which couples the incompressible Navier-Stokes equations and the Fokker-Planck equation. There is a lot of mathematical results about the FENE dumbbell model. M. Renardy \cite{Renardy} established the local well-posedness in Sobolev spaces with potential $\mathcal{U}(R)=(1-|R|^2)^{1-\sigma}$ for $\sigma>1$. Later, B. Jourdain, T. Leli\`{e}vre, and
C. Le Bris \cite{Jourdain} proved local existence of a stochastic differential equation with potential $\mathcal{U}(R)=-k\log(1-|R|^{2})$ in the case $k>3$ for a Couette flow. H. Zhang and P. Zhang \cite{Zhang-H} proved local well-posedness of (1.4) with $d=3$ in weighted Sobolev spaces. For the co-rotation case, F. Lin, P. Zhang, and Z. Zhang \cite{F.Lin} obtain a global existence results with $d=2$ and $k > 6$. If the initial data is perturbation around equilibrium, N. Masmoudi \cite{Masmoudi2008} proved global well-posedness of (1.4) for $k>0$. In the co-rotation case with $d=2$, he \cite{Masmoudi2008} obtained a global result for $k>0$ without any small conditions. In the co-rotation case, A. V. Busuioc, I. S. Ciuperca, D. Iftimie and L. I. Palade \cite{Busuioc} obtain a global existence result with only the small condition on $\psi_0$. The global existence of weak solutions in $L^2$ was proved recently by N. Masmoudi \cite{Masmoudi2013} under some entropy conditions. The $L^2$ decay of solutions was studied by M. Schonbek \cite{Schonbek} and she proved that the $L^2$ decay rate of the velocity is $(1+t)^{-\frac{d}{4}+\frac{1}{2}}$ with $d\geq 2$. More recently, W. Luo and Z. Yin \cite{Luo-Yin} improved Schonbek's result and showed that the decay rate is $(1+t)^{-\frac{d}{4}}$ with $d\geq 3$ and $\ln^{-l}(1+t)$ with $d=2$ for any $l\in\mathbb{N^+}$.

All the models we mentioned above are viscoelastic models. The main different between them is the constitutive equation. The constitutive equation of the Oldroyd-B model is very simple form $\sigma=G^TG$.  For the FENE dumbbell model, $\sigma$ is given by the Kramers expression i.e. $\sigma=\int_{B}R\otimes\mathcal{U}\psi dR$. In the Doi-Edwards polymer model, $\sigma$ can be deduced from a tensor $S$ denoting the orientation of the chains. The tensor $S$ is given by an integral constitutive law $S=-\displaystyle\int^{+\infty}_0\partial_{\tau}K(t,\tau,x,s)\mathscr{S}(G(t,\tau,x))d\tau$.

There are two time variables $t$ and $\tau$ in the Doi-Edwards polymer model. Usually, the unknown functions is depend on the current time $t$ and the past time $t'$. And $\tau=t-t'$ is the age of the polymer. In the paper \cite{Chupin2017}, the author consider the time variables $t$ and $\tau$ as two independent variables. In this paper, we still assume that $\tau$ is not depend on $t$.

Although there are a lot of research papers studied about the mathematical theory for viscoelastic models, the mathematical results for the Doi-Edwards polymer model remains very poor.  Recently, L. Chupin investigated \cite{Chupin2013,Chupin2014,Chupin2017} the global existence results of the Doi-Edwards polymer model. For a simplified model of (\ref{1.1}), L. Chupin \cite{Chupin2013,Chupin2014} proved the local existence and uniqueness of strong solutions with $x\in\mathbb{T}^d$ in Sobolev spaces. In \cite{Chupin2013}, L. Chupin shows that if the initial velocity $v_0$ and the viscosities ratio $w$ is small enough, then the local solutions is global in time. In \cite{Chupin2014}, L. Chupin obtained a global results without small conditions with $d=2$. In \cite{Chupin2017}, L. Chupin proved existence and uniqueness of local solutions of (\ref{1.1}) with $x\in\mathbb{T}^d$. When $d=2$, he \cite{Chupin2017} obtained a global result without small conditions.

In the system (\ref{1.1}), we observe that the boundary conditions $K(t,0, x,s)=1$ and $K(t,\tau,x,\frac{1}{2})=K(t,\tau,x,-\frac{1}{2})=0$ are not compatible. Hence, we consider the memory function $m=-\partial_\tau K$. On the other hand, we denote that $F=G-I$ and $\widetilde{\mathscr{S}}(F)=\mathscr{S}(G)=\mathscr{S}(F+I)$, and then (\ref{1.1}) can be written in the following form:
   \begin{align}\label{1.5}
\left\{
\begin{array}{ll}
Re [\partial_{t}v+(v\cdot\nabla)v]-(1-\omega)\Delta{v}+\nabla{P}=div~\sigma, ~~~~~~~div ~v=0,\\[1ex]
\sigma(t,x)=\omega\displaystyle\int^{\frac{1}{2}}_{-\frac{1}{2}}S(t,x,s)ds,~~S(t,x,s)=\displaystyle\int^{+\infty}_0m(t,\tau,x,s)\widetilde{\mathscr{S}}(F(t,\tau,x))d\tau,\\[1ex]
\partial_{t}F+(v\cdot\nabla)F+\frac{1}{We}\partial_\tau F=F\cdot\nabla v+\nabla v,  \\[1ex]
\partial_tm+(v\cdot\nabla)m+\frac{1}{We}\partial_\tau m+(\nabla v:\int^s_0S)\partial_sm-\frac{1}{We}\partial^2_sm=0, \\[1ex]
v(t,x)|_{t=0}=v_{0},\quad F(t,\tau,x)|_{t=0}=F_{0}(\tau,x) ,\quad m(t,\tau,x,s)|_{t=0}=m_0(\tau,x,s), \\[1ex]
m(t,\tau,x,s)|_{\tau=0}=0,\quad m(t,\tau,x,s)|_{s=\frac{1}{2}}=m(t,\tau,x,s)|_{s=-\frac{1}{2}}=0, \quad  F(t,\tau,x)|_{\tau=0}=0.\\[1ex]
\end{array}
\right.
\end{align}

In this paper, we are interested in the Cauchy problem of (\ref{1.5}) in the whole space $x\in\mathbb{R}^d$ with $d\geq 2$. To our best knowledge, there is no any well-posedness results for (\ref{1.5}) in the whole spaces $x\in\mathbb{R}^d$. We will use the Littlewood-Paley theory to study about the local well-posedness of (\ref{1.5}) in Besov spaces. This method allows us to assume that the initial data has lower regularity. By using a contraction argument, we also prove that if the initial velocity and the initial memory is small enough then the local solution exists global in time. 

\begin{rema}
In $\cite{Chemin2001}$, J. Chemin and N. Masmoudi's result implies that if the $L^\infty-$norm of $\tau$ is bounded and the velocity $v$ has finite energy, then the Navier-Stokes equations has a global solution with $d=2$. By virtue of the maximum principle, one can obtain the $L^\infty-$norm of $\tau$. Since $L^\infty\hookrightarrow L^2$ when $x\in\mathbb{T}^d$, it follows that the energy of velocity $v$ is finite. However, the embedding result $L^\infty\hookrightarrow L^2$ is not true for $x\in \mathbb{R}^d$. Hence, for the system (\ref{1.5}), there is no basic energy estimate when $x\in\mathbb{R}^d$. This is the reason that we can not obtain the global result without any small condition.
\end{rema}

Even if we assume that the initial velocity is small, there are still some difficult which is different with the other system. The main difficult is that the linear system of (\ref{1.5}) is quite different from the Oldroyd-B model. It seems that (\ref{1.5}) has no dissipative structure, so even though the initial data $(v_0,F_0,m_0)$ is small enough, one can't not obtain that the corresponding solution $(v,F,m)$ is small in all time. Therefore, the methods which used to deal with the Oldroyd-B model are invalid in this paper. However, we observe that the stress term is a integral form, this allows us to get some fine properties. Concretely, if we consider $m$ belongs to some suitable $\tau$-weighted function spaces, then the Besov norm of stress term will depend on this weighted. Under some condition on the weighted, we can see that the stress term will become small. Due to this observation, we only add the small condition on the initial velocity $v_0$. In order to apply the bootstrap argument to prove that the solution $v$ is also small for all time, we need to get a uniform estimate for $\tau$ independent of $t$.  This estimate is also a difficult, we must prove a maximum principle with $\tau$-weighted and the exponential decay of $m$. This result only holds true for the viscoelastic model with an integral constitutive law.

  The paper is organized as follows. In Section 2 we introduce some notations and our main results. In Section 3 we give some preliminaries which will be used in the sequel. In Section 4 we investigate the linear problem of (\ref{1.5}) and give some a priori estimates for solutions to (\ref{1.5}). In Section 5 we prove the local well-posedness of (\ref{1.5}) by using an approximate argument.  In Section 6 we prove the global well-posedness of (\ref{1.5}) by a contradiction argument.

\section{Notations and main results}
  In this section we introduce our main results and the notations that we shall use throughout the paper.

We use the symbol $L^p$ to denoted the $p$-integrable spaces. For a function $f$ depends only on the spatial variable $x$, we define that
$$\|f\|_{L^p_x}=(\int_{\mathbb{R}^d}|f|^pdx)^{\frac{1}{p}}.$$
For a function $g$ depends on $x$ and $s$, we define that
$$\|g\|_{L^p_{x,s}}=(\int^{\frac{1}{2}}_{-\frac{1}{2}}\int_{\mathbb{R}^d}|g|^pdxds)^{\frac{1}{p}}.$$
For simplicity, sometimes we will drop the subscript and use the symbol $L^p$ to denote $L^p_x$ and $L^p_{x,s}$ if there is no ambiguity.

The symbol $\mathscr{F}$ is the Fourier transform and $\mathscr{F}^{-1}$ represents its inverse.  Moreover, we use $S$ to denote the Schwartz space $S(\mathbb{R}^d)$ and $S'$ is the dual space of $S$. Finally, let us define a subspace of $S'$, namely $S'_h$:
$$S'_h=\big\{u\in S'\big| \lim_{\lambda\rightarrow\infty}\|\theta(\lambda D)u\|_{L^\infty}=0,~~\forall \theta\in C^\infty_0\big\},$$
where the pseudo-differential operator $\theta(D) $ is defined by $\theta(D) a(x)\triangleq \mathscr{F}^{-1}(\theta(\xi) \mathscr{F} a(\xi))$.\\
  Next we introduce the Littlewood-Paley decomposition and homogeneous Besov spaces (see \cite{Bahouri2011} for more details).

Let $\mathcal{C}$ be the annulus $\{\xi\in\mathbb{R}^{d}\big|\frac{3}{4}\leq|\xi|\leq\frac{8}{3}\}.$ There exists radial function $\chi$ and $\varphi$, valued in the interval $[0,1]$, such that
\begin{align}
\forall \xi\in\mathbb{R}^d,  \chi(\xi)+\sum_{j\geq 0}\varphi(2^{-j}\xi)=1,
\end{align}
\begin{align}
\forall\xi\in\mathbb{R}^{d}\backslash\{0\},~\sum_{j\in\mathbb{Z}}\varphi(2^{-j}\xi)=1,
\end{align}
\begin{align}
|j-j'|\geq2\Rightarrow Supp ~\varphi(2^{-j}\xi)\cap Supp ~\varphi(2^{-j'}\xi)=\emptyset.
\end{align}
 The homogeneous dyadic blocks $\dot{\Delta}_{j}$ are defined by
\begin{align} \dot{\Delta}_{j}u=\varphi(2^{-j}D)u=2^{jd}\int_{\mathbb{R}^{d}}h(2^{j}y)u(x-y)dy,
\end{align}
\begin{align}
\dot{S}_{j}u=\chi(2^{-j}D)u=\int_{\mathbb{R}^{d}}\widetilde{h}(2^{j}y)u(x-y)dy,
\end{align}
where $h=\mathscr{F}^{-1}\varphi$ and $\widetilde{h}=\mathscr{F}^{-1}\chi$.

   The homogeneous Besov space is denoted by $\dot{B}^{s}_{p,r}$ i.e.,
$$\dot{B}^{s}_{p,r}=\big\{u\in S'_{h}\big{|}\|u\|_{\dot{B}^{s}_{p,r}}=\|2^{js}\|\dot{\Delta}_{j}u\|_{L^{p}_x}\|_{l^r}<\infty\big\},$$

  The notation $C([0,T];X)$ denotes the space of continuous functions on $[0,T]$ with values in the Banach space $X$.  And $L^{\rho}(0,T;X)$ denotes the space of $\rho$-integrable functions on $(0,T)$ with values in the Banach space $X$.

  Note that $\tau\in[0,+\infty)$. We denoted the space $L^\infty_{\tau}(X)$ endowed with the norm
  $$\|f\|_{L^\infty_\tau(X)}=ess~\sup_{\tau\in\mathbb{R}^+}\|f(\tau,\cdot)\|_{X}.$$

  Also we will use the following spaces
$$\widetilde{L}^{\rho}_{t}(\dot{B}^{s}_{p,r})=\big\{u\in S'\big{|}\|u\|_{\widetilde{L}^{\rho}_{t}(\dot{B}^{s}_{p,r})}=\|2^{js}\|\dot{\Delta_{j}}u\|_{L^{\rho}_{t}(L^{p})}\|_{l^r}<\infty\big\},$$

$$\dot{\mathcal{B}}^{s}_{p,r}=\big\{u\in S'_{h}\big{|}\|u\|_{\dot{B}^{s}_{p,r}}=\|2^{js}\|\dot{\Delta}_{j}u\|_{L^{p}_{x,s}}\|_{l^r}<\infty\big\}. $$
By virtue of the Minkowski inequality, we get
$$\widetilde{L}^2_T(\dot{B}^{s}_{p,1})\hookrightarrow L^2_T(\dot{B}^{s}_{p,1})\hookrightarrow L^1_T(\dot{B}^{s}_{p,1}).$$

For simplicity, we assume that the Reynolds number $Re=1$ and the Weissenberg number $We=1$. The notation $a\lm b$ means that there is a universal constant $C$ such that $a\leq Cb$.

Throughout this paper, we assume that there exists a constant $\gamma>0$ such that
$$det~G_0=det~(F_0+I)\geq \gamma>0.$$
The above assumption will leads to a lower bound for $G=F+I$ (See Lemma 1 in \cite{Chupin2017}). The lower bound for $F-I$ is crucial to estimate the Besov norm of the function $\widetilde{\mathscr{S}}(F(t,\tau,x))$ (See Proposition \ref{P4} in Chapter 4).

Now we state our main results as follows.
\begin{theo}
Let $2\leq p<2d$ and $\lambda>0$. Assume that $v_{0}\in \dot{B}^{\frac{d}{p}}_{p,1}\cap \dot{B}^{-1+\frac{d}{p}}_{p,1},~F_{0}\in L^\infty_\tau(\dot{B}^{\frac{d}{p}}_{p,1}\cap \dot{B}^{-1+\frac{d}{p}}_{p,1}),$ and $(\lambda+\tau)^2m_0\in L^\infty_\tau(\dot{\mathcal{B}}^{\frac{d}{p}}_{p,1}\cap \dot{B}^{-1+\frac{d}{p}}_{p,1})$. Then there exist some $T^{*}>0$ and a unique solution $(v,F,m,P)$ of $(1.5)$ such that
$$v\in C([0,T^{*});\dot{B}^{\frac{d}{p}}_{p,1}\cap\dot{B}^{-1+\frac{d}{p}}_{p,1})\cap \widetilde{L^2}(0,T^{*};\dot{B}^{1+\frac{d}{p}}_{p,1}),~~~\nabla P\in L^\infty([0,T^{*});\dot{B}^{\frac{d}{p}}_{p,1}\cap\dot{B}^{-1+\frac{d}{p}}_{p,1}),$$
$$ F\in C([0,T^{*});L^\infty_\tau(\dot{B}^{\frac{d}{p}}_{p,1}\dot{B}^{-1+\frac{d}{p}}_{p,1})),~~~(\lambda+\tau)^2m\in C([0,T^{*});L^\infty_\tau(\dot{\mathcal{B}}^{\frac{d}{p}}_{p,1}\cap\dot{B}^{-1+\frac{d}{p}}_{p,1})).$$
\end{theo}

\begin{theo}
Under the assumption of Theorem 2.1. If there exists a constant $\ep$ such that $$\|v_{0}\|_{\dot{B}^{-1+\frac{d}{p}}_{p,1}\cap \dot{B}^{\frac{d}{p}}_{p,1}}\leq \ep,~~~\frac{1}{\lambda}\leq \ep^2,  $$
then the solution constructed in Theorem 2.1 is global.
\end{theo}
\begin{rema}
If we replace the small condition $\frac{1}{\lambda}\leq \ep^2$ by the viscosities ratio $w\leq \ep^2$, then the global result still holds true. The proof is similar to that of Theorem 2.2. 
\end{rema}

\section{Preliminaries}
In this section we introduce some useful lemmas which will be used in the sequel. For more details, one can refer to Section 2 in \cite{Bahouri2011}.
\subsection{Propositions of Besov spaces}
Firstly we introduce the Bernstein inequalities.
\begin{lemm}\label{L1}
\cite{Bahouri2011} Let $\mathcal{C}$ be an annulus and $B$ a ball. For any nonnegative integer $k$, any couple $(p, q)$ in $[1,\infty]^{2}$ with $q \geq p\geq 1$, and
any function $u$ of $L^{p}$, we have
$$Supp ~\widehat{u} \subseteq \lambda B \Rightarrow \|D^{k}u\|_{L^{q}}\triangleq \sup_{|\alpha|\leq k} \|\partial ^{\alpha}u\|_{L^{q}}\lm \lambda^{k+d(\frac{1}{p}-\frac{1}{q})}\|u\|_{L^{p}},$$
$$Supp ~\widehat{u} \subseteq \lambda\mathcal{C} \Rightarrow \lambda^{k}\|u\|_{L^{p}} \lm \|D^{k}u\|_{L^{p}}\lm \lambda^{k}\|u\|_{L^{p}}.$$
$$Supp ~\widehat{u} \subseteq \lambda\mathcal{C} \Rightarrow  \|e^{t\triangle}u\|_{L^{p}}\lm e^{-ct\lambda^{2}}\|u\|_{L^{p}}.$$
\end{lemm}

The following proposition is about the embedding for Besov spaces.
\begin{prop}\label{P1}
\cite{Bahouri2011} Let $1\leq p_{1} \leq p_{2} \leq \infty$ and $1\leq r_{1} \leq r_{2} \leq \infty$, and let $s$ be a real number. Then we have
$$\dot{B}^{s}_{p_{1},r_{1}}\hookrightarrow \dot{B}^{s-d(\frac{1}{p_{1}}-\frac{1}{p_{2}})}_{p_{2},r_{2}}, \quad \dot{B}^{\frac{d}{p}}_{p,1}\hookrightarrow L^\infty.$$
\end{prop}

Next we introduce the Bony decomposition.
\begin{defi}\label{D1}
\cite{Bahouri2011} For functions $u$ and $v$, Bony's decomposition in the homogeneous context is defined by
$$uv=\dot{T}_{u}v+\dot{R}(u,v)+\dot{T}_{v}u,$$
where
$$\dot{T}_{u}v\triangleq\sum_{j}\dot{S}_{j-1}u\dot{\Delta}_{j}v,\quad \dot{R}(u,v)\triangleq\sum_{|k-j|\leq1}\dot{\Delta}_{k}u\dot{\Delta}_{j}v.$$
\end{defi}

\begin{prop}\label{P2}
\cite{Bahouri2011} Let $s\in\mathbb{R}$, $t<0$ and $(p, r_{1}, r_{2})$
in $[1,\infty]^{3}$, then
$$\|\dot{T}_{u}v\|_{\dot{B}^{s}_{p,r}}\lm \|u\|_{L^{\infty}}\|v\|_{\dot{B}^{s}_{p,r}},\quad  \|\dot{T}_{u}v\|_{\dot{B}^{s+t}_{p,r}}\lm \|u\|_{\dot{B}^{t}_{p,r_{1}}}\|v\|_{\dot{B}^{s}_{p,r_{2}}},$$
where $\frac{1}{r}=\min\{1,\frac{1}{r_{1}}+\frac{1}{r_{2}}\}.$
\end{prop}

\begin{prop}\label{P3}
\cite{Bahouri2011} Let $(s_{1}, s_{2})$ be in $\mathbb{R}^2$ and $(p_{1}, p_{2}, r_{1}, r_{2})$ be in $[1,\infty]^{4}$. Assume that
$$\frac{1}{p}=\frac{1}{p_{1}}+\frac{1}{p_{2}}~~and~~\frac{1}{r}=\frac{1}{r_{1}}+\frac{1}{r_{2}} .$$
If $s_{1}+s_{2}>0$, there exists a constant $C$ such that:
$$\|\dot{R}(u,v)\|_{\dot{B}^{s_{1}+s_{2}}_{p,r}}\leq \frac{C^{s_{1}+s_{2}+1}}{s_{1}+s_{2}}\|u\|_{\dot{B}^{s_{1}}_{p_{1},r_{1}}}\|v\|_{\dot{B}^{s_{2}}_{p_{2},r_{2}}}.$$
If $r=1$ and $s_{1}+s_{2}=0$, then we have,
$$\|\dot{R}(u,v)\|_{\dot{B}^{0}_{p,\infty}}\leq C\|u\|_{\dot{B}^{s_{1}}_{p_{1},r_{1}}}\|v\|_{\dot{B}^{s_{2}}_{p_{2},r_{2}}}.$$
\end{prop}

\begin{coro}\label{C1}
\cite{Bahouri2011} If $(s, p, r)$ in $(0,\infty]\times[1,\infty]^{2}$, the
space $L^{\infty}\cap \dot{B}^{s}_{p,r}$ is an algebra, moreover we have
$$\|uv\|_{\dot{B}^{s}_{p,r}}\lm \|u\|_{L^{\infty}}\|v\|_{\dot{B}^{s}_{p,r}}+\|u\|_{\dot{B}^{s}_{p,r}}\|v\|_{L^{\infty}}.$$
If $s=\frac{d}{p},~r=1$, we have
$$\|uv\|_{\dot{B}^{\frac{d}{p}}_{p,1}}\lm \|u\|_{\dot{B}^{\frac{d}{p}}_{p,1}}\|v\|_{\dot{B}^{\frac{d}{p}}_{p,1}}.$$
\end{coro}

\begin{lemm}\label{L2}
\cite{Bahouri2011} Let $p\in[1,+\infty]$, $-1<\sigma\leq \frac{d}{p}$ and $v$ be a vector field. For any function $f\in \dot{B}^{\sigma}_{p,1}$, we have
$$\|2^{j\sigma}\|R_j\|_{L^p}\|_{l^1}\lm \|v\|_{\dot{B}^{1+\frac{d}{p}}_{p,1}}\|f\|_{\dot{B}^{\sigma}_{p,1}},$$
where $R_j=[v\nabla,\dot{\Delta_j}]f$.
\end{lemm}

The following product laws will be used constantly.
\begin{lemm}\label{L3}
\cite{Paicu-Zhang2012}Let $p\geq 1$, and $s_1\leq \frac{d}{p}, s_2\leq \frac{d}{p}$ with $s_1+s_2>d\max(0,\frac{2}{p}-1)$. If $a\in\dot{B}^{s_1}_{p,1}, \quad b\in\dot{B}^{s_2}_{p,1}$, then $ab\in\dot{B}^{s_1+s_2-\frac{d}{p}}_{p,1}$, and
$$\|ab\|_{\dot{B}^{s_1+s_2-\frac{d}{p}}_{p,1}}\lm \|a\|_{\dot{B}^{s_1}_{p,1}}\|b\|_{\dot{B}^{s_2}_{p,1}}.$$
\end{lemm}

\begin{lemm}\label{L4}
\cite{Bahouri2011} Let $f$ be a smooth function on $\mathbb{R}$ which vanishes at $0$. Let $s$ be a real number and $(p,r)\in[1,+\infty]^2$ satisfies that
\[s<\frac{d}{p}, \quad or \quad s=\frac{d}{p},~r=1.\]
For any function $u\in \dot{B}^s_{p,r}\cap L^\infty$, we have
$$\|f\circ u\|_{\dot{B}^s_{p,r}}\leq C(f',\|u\|_{L^\infty})\|u\|_{\dot{B}^s_{p,r}}.$$
\end{lemm}
\subsection{A Gronwall lemma with two variables}
The equations about the functions $m$ and $F$ are involve two independent time variables $t$ and $\tau$. In order to obtain some a priori estimates, we need the following Gronwall lemma with two variables.
\begin{lemm}\label{L5}
\cite{Chupin2013}Let $f,g:\mathbb{R}^{+}\mapsto \mathbb{R}^{+}$ be two positive and locally integrable functions. If a function $y:\mathbb{R}^{+}\times\mathbb{R}^{+}\mapsto \mathbb{R}^{+}$ satisfies that
$$\partial_ty(t,\tau)+\partial_\tau y(t,\tau)\leq f(t)y(t,\tau)+g(t), \quad \forall (t,\tau)\in(0,T)\times\mathbb{R}^{+},$$
then we have
$$y(t,\tau)\leq \bigg[\zeta(t,\tau)+\int^t_0g(t')dt'\bigg]\exp\bigg(\int^t_0f(t')dt'\bigg), \quad \forall (t,\tau)\in(0,T)\times\mathbb{R}^{+},$$
where $\zeta(t,\tau)=\begin{cases}
y(\tau-t,0)\quad \text{if} \quad t\leq \tau,\\
y(0,t-\tau)\quad \text{if}  \quad t>\tau.
\end{cases}$
\end{lemm}

\subsection{A maximum principle}
The following maximum principle will be used get the bound for the stress tensor.
\begin{lemm}\label{L6}
\cite{Chupin2017}If $\partial_sg\in L^1(0,T;L^\infty_{x,s})$ and $v$ is divergence free, then the solution $f(t,\tau,x,s)$ to the following system:
\begin{equation}
\left\{
\begin{array}{ll}
\partial_tf+\partial_\tau f+(v\cdot\nabla) f+g\partial_sf-\partial^2_sf=0, \\
f|_{t=0}=f_0, \quad f|_{\tau=0}=f_1, \quad f|_{s=\frac{1}{2}}=f|_{s=-\frac{1}{2}}=0,
\end{array}
\right.
\end{equation}
satisfies the following maximum principle on $(0,T)\times\mathbb{R}^{+}\times\mathbb{R}^d\times(-\frac{1}{2},\frac{1}{2})$:
$$\min\{\inf_{\tau,x,s} f_0,~\inf_{\tau,x,s} f_1\}\leq f\leq \max\{\sup_{\tau,x,s} f_0,~\sup_{\tau,x,s} f_1\}.$$
\end{lemm}

\begin{lemm}\label{L7}
Suppose that $\lambda>\frac{1}{4\pi^2}$. If $\partial_sg\in L^1(0,T;L^\infty_{x,s})$ and $v$ is divergence free, then the solution $f(t,\tau,x,s)$ to the following system:
\begin{equation}
\left\{
\begin{array}{ll}
\partial_tf+\partial_\tau f+(v\cdot\nabla) f+g\partial_sf-\frac{1}{\lambda+\tau}f-\partial^2_sf=0, \\
f|_{t=0}=f_0, \quad f|_{\tau=0}=f_1, \quad f|_{s=\frac{1}{2}}=f|_{s=-\frac{1}{2}}=0,
\end{array}
\right.
\end{equation}
satisfies the following maximum principle on $(0,T)\times\mathbb{R}^{+}\times\mathbb{R}^d\times(-\frac{1}{2},\frac{1}{2})$:
$$\min\{\inf_{\tau,x,s} f_0,~\inf_{\tau,x,s} f_1\}\leq f\leq \max\{\sup_{\tau,x,s} f_0,~\sup_{\tau,x,s} f_1\}.$$
\begin{proof}
Note that $4\pi^2$ is the first eigenvalue of the operator $-\partial^2_s$ associated with the above boundary conditions. Since $\frac{1}{\lambda+\tau}<\frac{1}{\lambda}<4\pi^2$, then one can obtain the maximum principle by virtue of the similar proof in \cite{Chupin2017}.
\end{proof}
\end{lemm}

\subsection{The estimate for stress tensor}
In this section, we will give the bound for the stress tensor.  First we bound the function $\mathscr{S}$.
\begin{prop}\label{P4}
The function $\mathscr{S}$ is defined by the relation (1.2) and $\widetilde{\mathscr{S}}(F)=\mathscr{S}(G)=\mathscr{S}(F+I)$. Let $s$ be a real number and $(p,r)\in[1,+\infty]^2$ satisfies that
\[s<\frac{d}{p}, \quad or \quad s=\frac{d}{p},~r=1.\]
The following properties hold true:
$$|\widetilde{\mathscr{S}}(F)|\leq C,\quad |\widetilde{\mathscr{S}}'(F)|\leq \frac{C}{|F+I|}\quad \|\widetilde{\mathscr{S}}(F)\|_{\dot{B}^{s}_{p,r}}\leq C(\widetilde{|\mathscr{S}}'|,\|F\|_{L\infty})\|F\|_{\dot{B}^{s}_{p,r}}.$$
\begin{proof}
The estimates for $|\widetilde{\mathscr{S}}(G)|$ and $|\widetilde{\mathscr{S}}'(G)|$ can be found in Proposition 1 in \cite{Chupin2017}. By virtue of the relation (1.2), one can see that $\widetilde{\mathscr{S}}(0)=\mathscr{S}(I)=0$. Thanks to Lemma \ref{L4}, we complete the proof.
\end{proof}
\end{prop}

\begin{rema}
In Chapter 2, we have assumed that $det(F_0+I)\geq \gamma>0$ which leads to $|F+I|\geq C_\gamma>0$(See Lemma 1 in \cite{Chupin2017}). By virtue of the above proposition, we deduced that $|\widetilde{\mathscr{S}}'(F)|\leq \frac{C}{|F+I|}\leq \frac{C}{C_\gamma}$. Thus we can bounds for the Besov norm of $\widetilde{\mathscr{S}}(F)$ if the Besov norm of $F$ is finite.
\end{rema}

\begin{coro}\label{C2}
Let $(\lambda+\tau)^2m\in L^\infty_\tau(\dot{B}^{\frac{d}{p}}_{p,1})$ and $F\in L^\infty_\tau(\dot{B}^{\frac{d}{p}}_{p,1})$ with $p\in[1,+\infty]$. Then
\[\|S\|_{\dot{B}^{\frac{d}{p}}_{p,1}}\lm \frac{1}{\lambda}\|(\lambda+\tau)^2m\|_{L^\infty_\tau(\dot{B}^{\frac{d}{p}}_{p,1})}\|F\|_{L^\infty_\tau(\dot{B}^{\frac{d}{p}}_{p,1})},\]
\[\|S\|_{L^\infty_{x,s}}\lm \frac{1}{\lambda}\|(\lambda+\tau)^2m\|_{L^\infty_\tau(L^\infty_{x,s})},\]
where $S=\displaystyle\int^{+\infty}_0m\widetilde{\mathscr{S}}(F)d\tau$.
\begin{proof}
By virtue of Corollary \ref{C1} and Proposition \ref{P4}, we obtain
\begin{align*}
\|S\|_{\dot{B}^{\frac{d}{p}}_{p,1}}&=\|\displaystyle\int^{+\infty}_0m\widetilde{\mathscr{S}}(F)d\tau\|_{\dot{B}^{\frac{d}{p}}_{p,1}}\leq\displaystyle\int^{+\infty}_0(\lambda+\tau)^{-2}\|(\lambda+\tau)^2m\widetilde{\mathscr{S}}(F)\|_{\dot{B}^{\frac{d}{p}}_{p,1}}d\tau\\
&\lm \|(\lambda+\tau)^2m\|_{L^\infty_\tau(\dot{B}^{\frac{d}{p}}_{p,1})}\|F\|_{L^\infty_\tau(\dot{B}^{\frac{d}{p}}_{p,1})}\int^\infty_0(\lambda+\tau)^{-2}d\tau=\frac{1}{\lambda}\|(\lambda+\tau)^2m\|_{L^\infty_\tau(\dot{B}^{\frac{d}{p}}_{p,1})}\|F\|_{L^\infty_\tau(\dot{B}^{\frac{d}{p}}_{p,1})},
\end{align*}
and
\begin{align*}
\|S\|_{L^\infty_{x,s}}&=\|\displaystyle\int^{+\infty}_0m\widetilde{\mathscr{S}}(F)d\tau\|_{L^\infty}\leq\displaystyle\int^{+\infty}_0(\lambda+\tau)^{-2}\|(\lambda+\tau)^2m\widetilde{\mathscr{S}}(F)\|_{L^\infty}d\tau\\
&\leq \|(\lambda+\tau)^2m\|_{L^\infty_\tau(L^\infty_{x,s})}\|\widetilde{\mathscr{S}}(F)\|_{L^\infty}\int^\infty_0(\lambda+\tau)^{-2}d\tau \lm \frac{1}{\lambda}\|(\lambda+\tau)^2m\|_{L^\infty_\tau(L^\infty_{x,s})}.
\end{align*}

\end{proof}
\end{coro}

\section{Linear problem and a priori estimates}
In this section, we assume that the vector field $u$ and the tensor $S$ are given, and then consider the following linearized equations for (\ref{1.5}):
   \begin{align}\label{4.1}
\left\{
\begin{array}{ll}
\partial_{t}v+(u\cdot\nabla)v-(1-\omega)\Delta{v}+\nabla{P}=div~\sigma, ~~~~~~~div ~v=0,\\[1ex]
\sigma(t,x)=\omega\displaystyle\int^{\frac{1}{2}}_{-\frac{1}{2}}S(t,x,s)ds, \\[1ex]
\partial_{t}F+(u\cdot\nabla)F+\partial_\tau F=F\cdot\nabla u+\nabla u,  \\[1ex]
\partial_tm+(u\cdot\nabla)m+\partial_\tau m+(\nabla u:\int^s_0S)\partial_sm-\partial^2_sm=0, \\[1ex]
v(t,x)|_{t=0}=v_{0},\quad F(t,\tau,x)|_{t=0}=F_{0}(\tau,x) ,\quad m(t,\tau,x,s)|_{t=0}=m_0(\tau,x,s), \\[1ex]
m(t,\tau,x,s)|_{\tau=0}=0,\quad m(t,\tau,x,s)|_{s=\frac{1}{2}}=m(t,\tau,x,s)|_{s=-\frac{1}{2}}=0, \quad  F(t,\tau,x)|_{\tau=0}=0.\\[1ex]
\end{array}
\right.
\end{align}

\subsection{Solutions to the linear equations for $F$}
\begin{prop}\label{P4.1}
Let $p\in[1,+\infty)$. Assume that $u\in \widetilde {L^2}(0,T;\dot{B}^{1+\frac{d}{p}}_{p,1})\hookrightarrow L^1(0,T;\dot{B}^{1+\frac{d}{p}}_{p,1})$ with $div~u=0$ and $F_0\in L^\infty_\tau(\dot{B}^{\frac{d}{p}}_{p,1})$. Then the following equation
\begin{align}\label{4.2}
\left\{
\begin{array}{ll}
\partial_{t}F+(u\cdot\nabla)F+\partial_\tau F=F\cdot\nabla u+\nabla u,  \\[1ex]
F(t,\tau,x)|_{t=0}=F_{0}(\tau,x), \quad F(t,\tau,x)|_{\tau=0}=0,
\end{array}
\right.
\end{align}
has a unique solution $F$ in $C([0,T];L^\infty_\tau(\dot{B}^{\frac{d}{p}}_{p,1}))$. Moreover, we have
\begin{align}\label{4.3}
\sup_{t\in[0,T]}\|F\|_{L^\infty_\tau(\dot{B}^{\frac{d}{p}}_{p,1})}\lm \exp(\int^T_0\|u\|_{\dot{B}^{1+\frac{d}{p}}_{p,1}}dt')(\|F_0\|_{L^\infty_\tau(\dot{B}^{\frac{d}{p}}_{p,1})}+\int^T_0\|u\|_{\dot{B}^{1+\frac{d}{p}}_{p,1}}dt').
\end{align}
\begin{proof}
Firstly, we prove the a priori estimate (\ref{4.3}) for $F$. Applying $\dot{\Delta}_j$ to the both sides of \ref{4.2}, we obtain that
\begin{align}
\partial_{t}\dot{\Delta}_jF+(u\cdot\nabla)\dot{\Delta}_jF+\partial_\tau \dot{\Delta}_jF=R_j+\dot{\Delta}_j(F\cdot\nabla u)+\dot{\Delta}_j\nabla u,
\end{align}
where $R_j=[u\nabla, \dot{\Delta_j}]F$. By virtue of the standard $L^p$ estimate, we see that
\begin{align*}
\partial_t\|\dot{\Delta_j}F\|_{L^p}+\partial_\tau\|\dot{\Delta_j}F\|_{L^p}\leq \|R_j\|_{L^p}+\|\dot{\Delta}_j(F\cdot\nabla u)\|_{L^p}+\|\dot{\Delta}_j\nabla u\|_{L^p}.
\end{align*}
Multiplying $2^{j\frac{d}{p}}$ and taking the $l^1$-norm in the both sides of the above inequality, we deduce that
\begin{align*}
\partial_t\|F\|_{\dot{B}^{\frac{d}{p}}_{p,1}}+\partial_\tau\|F\|_{\dot{B}^{\frac{d}{p}}_{p,1}}\leq \|2^{j\frac{d}{p}}\|R_j\|_{L^p}\|_{l^1}+\|F\cdot\nabla u\|_{\dot{B}^{\frac{d}{p}}_{p,1}}+\|\nabla u\|_{\dot{B}^{\frac{d}{p}}_{p,1}}.
\end{align*}
According to Corollary \ref{C1} and Lemma \ref{L2}, we get
\begin{align*}
\partial_t\|F\|_{\dot{B}^{\frac{d}{p}}_{p,1}}+\partial_\tau\|F\|_{\dot{B}^{\frac{d}{p}}_{p,1}}\lm \|F\|_{\dot{B}^{\frac{d}{p}}_{p,1}}\|\nabla u\|_{\dot{B}^{\frac{d}{p}}_{p,1}}+\|\nabla u\|_{\dot{B}^{\frac{d}{p}}_{p,1}}.
\end{align*}
Thanks to Lemma \ref{L5}, we obtain that
\begin{align*}
\sup_{t\in[0,T]}\|F\|_{L^\infty_\tau(\dot{B}^{\frac{d}{p}}_{p,1})}\lm \exp(\int^T_0\|u\|_{\dot{B}^{1+\frac{d}{p}}_{p,1}}dt')(\|F_0\|_{L^\infty_\tau(\dot{B}^{\frac{d}{p}}_{p,1})}+\int^T_0\|u\|_{\dot{B}^{1+\frac{d}{p}}_{p,1}}dt').
\end{align*}
Now, we prove the existence of solution to (\ref{4.2}). Assume that $(F^n_0, u^n)$ are smooth approximation of $(F_0, u)$ such that $F^n\rightarrow F_0$ in $L^\infty_\tau(\dot{B}^{\frac{d}{p}}_{p,1})$ and $u^n\rightarrow u$ in $ L^1(0,T;\dot{B}^{1+\frac{d}{p}}_{p,1})$. By virtue of the transport equation theory, for any fixed $\tau$ and $n$, there exists a function $F^n(t,x,\tau)$ satisfies the following equations:
\begin{align}\label{4.5}
\left\{
\begin{array}{ll}
\partial_{t}F^n+(u^n\cdot\nabla)F^n+\partial_\tau F^n=F^n\cdot\nabla u^n+\nabla u^n,  \\[1ex]
F|_{t=0}=F^n_{0}(\tau,x).
\end{array}
\right.
\end{align}
By virtue of (\ref{4.3}), we deduce that $F^n$ is uniformly bounded in $L^\infty(0,T;L^\infty_\tau(\dot{B}^{\frac{d}{p}}_{p,1}))$. The Fatou property for Besov spaces implies that there exists a subsequence $F^n_k$ weakly converges to some $F$ in  $L^\infty(0,T;L^\infty_\tau(\dot{B}^{\frac{d}{p}}_{p,1}))$. Since $\dot{B}^{\frac{d}{p}}_{p,1}$ is an algebra, taking the limit in (\ref{4.5}) as $n\rightarrow\infty$, it follows that $F$ is the solution to (\ref{4.2}). Moreover, one can check that $F\in C([0,T];L^\infty_\tau(\dot{B}^{\frac{d}{p}}_{p,1}))$. The method is similar as in Theorem 3.19 of \cite{Bahouri2011}, and we omit the details here. Using the fact that (\ref{4.2}) is a linear equation, the uniqueness can be proved by the estimate (\ref{4.3}).
\end{proof}
\end{prop}

\begin{prop}\label{P4.2}
Let ${2\leq p< 2d}$. Assume that $u\in L^1(0,T;\dot{B}^{1+\frac{d}{p}}_{p,1})$ with $div~u=0$, $w\in L^1(0,T;\dot{B}^{\frac{d}{p}}_{p,1})$, $H\in L^\infty(0,T;\dot{B}^{\frac{d}{p}}_{p,1})$ and $F_0\in L^\infty_\tau(\dot{B}^{-1+\frac{d}{p}}_{p,1})$. $F\in L^\infty(0,T;L^\infty_\tau(\dot{B}^{-1+\frac{d}{p}}_{p,1}))$ is the solution of the following equation
\begin{align*}
\left\{
\begin{array}{ll}
\partial_{t}F+(u\cdot\nabla)F+\partial_\tau F=w\cdot\nabla H+F\cdot \nabla u +H\cdot \nabla w+\nabla w,  \\[1ex]
F|_{t=0}=F_{0}(\tau,x), \quad F(t,\tau,x)|_{\tau=0}=0.
\end{array}
\right.
\end{align*}
Then, we have
\begin{align*}
\sup_{t\in[0,T]}\|F\|_{L^\infty_\tau(\dot{B}^{-1+\frac{d}{p}}_{p,1})}\lm \exp(\int^T_0\|u\|_{\dot{B}^{1+\frac{d}{p}}_{p,1}}dt')\bigg[\|F_0\|_{L^\infty_\tau(\dot{B}^{-1+\frac{d}{p}}_{p,1})}+\int^T_0\|w\|_{\dot{B}^{\frac{d}{p}}_{p,1}}(1+\|H\|_{\dot{B}^{\frac{d}{p}}_{p,1}})dt'\bigg].
\end{align*}
\begin{proof}
By the similar argument as in the Proposition \ref{P4.1}, we obtain that
\begin{align*}
\partial_t\|F\|_{\dot{B}^{-1+\frac{d}{p}}_{p,1}}+\partial_\tau\|F\|_{\dot{B}^{-1+\frac{d}{p}}_{p,1}}&\leq \|2^{j(-1+\frac{d}{p})}\|R_j\|_{L^p}\|_{l^1}+\|w\cdot\nabla H\|_{\dot{B}^{-1+\frac{d}{p}}_{p,1}}\\&
+\|H\cdot\nabla w\|_{\dot{B}^{-1+\frac{d}{p}}_{p,1}}+\|F\cdot\nabla u\|_{\dot{B}^{-1+\frac{d}{p}}_{p,1}}+\|\nabla w\|_{\dot{B}^{-1+\frac{d}{p}}_{p,1}}.
\end{align*}
Note that $2\leq p<2d$. Using Lemma \ref{L3} with $s_1=-1+\frac{d}{p}$ and $s_2=\frac{d}{p}$, we deduce that
\[\|w\cdot\nabla H\|_{\dot{B}^{-1+\frac{d}{p}}_{p,1}}\lm \|w\|_{\dot{B}^{\frac{d}{p}}_{p,1}}\|\nabla H\|_{\dot{B}^{-1+\frac{d}{p}}_{p,1}}\lm \|w\|_{\dot{B}^{\frac{d}{p}}_{p,1}}\|H\|_{\dot{B}^{\frac{d}{p}}_{p,1}},\]
\[\|H\nabla w\|_{\dot{B}^{-1+\frac{d}{p}}_{p,1}}\lm \|H\|_{\dot{B}^{\frac{d}{p}}_{p,1}}\|\nabla w\|_{\dot{B}^{-1+\frac{d}{p}}_{p,1}}\lm \|w\|_{\dot{B}^{\frac{d}{p}}_{p,1}}\|H\|_{\dot{B}^{\frac{d}{p}}_{p,1}}, \]
\[\|F\cdot\nabla u\|_{\dot{B}^{-1+\frac{d}{p}}_{p,1}}\lm \|F\|_{\dot{B}^{-1+\frac{d}{p}}_{p,1}}\|\nabla u\|_{\dot{B}^{\frac{d}{p}}_{p,1}}\lm \|u\|_{\dot{B}^{1+\frac{d}{p}}_{p,1}}\|F\|_{\dot{B}^{-1+\frac{d}{p}}_{p,1}},\]
which along with Lemma \ref{L2} yield that
  \begin{align*}
\partial_t\|F\|_{\dot{B}^{-1+\frac{d}{p}}_{p,1}}+\partial_\tau\|F\|_{\dot{B}^{-1+\frac{d}{p}}_{p,1}}&\lm \|u\|_{\dot{B}^{1+\frac{d}{p}}_{p,1}}\|F\|_{\dot{B}^{-1+\frac{d}{p}}_{p,1}}+\|w\|_{\dot{B}^{\frac{d}{p}}_{p,1}}(1+\|H\|_{\dot{B}^{\frac{d}{p}}_{p,1}}).
\end{align*}
Taking advantage of Lemma \ref{L5}, we obtain the desired result.
\end{proof}
\end{prop}

\subsection{Solutions to the linear equations for $m$}
\begin{prop}\label{P4.3}
Let $2\leq p<+\infty$. Assume that $u\in \widetilde{L^2}(0,T;\dot{B}^{1+\frac{d}{p}}_{p,1})$ with $div~u=0$, $S\in L^\infty(0,T;\dot{\mathcal{B}}^{\frac{d}{p}}_{p,1})$, $(\lambda+\tau)^2m_0\in L^\infty_\tau(\dot{\mathcal{B}}^{\frac{d}{p}}_{p,1})$. Then the following equation 
\begin{align}\label{4.6}
\left\{
\begin{array}{ll}
\partial_tm+(u\cdot\nabla)m+\partial_\tau m+(\nabla u:\int^s_0S)\partial_sm-\partial^2_sm=0,,  \\[1ex]
m|_{t=0}=m_{0}(\tau,x,s), \quad m(t,\tau,x,s)|_{\tau=0}=0, \quad m(t,\tau,x,s)|_{s=\frac{1}{2}}=m(t,\tau,x,s)|_{s=-\frac{1}{2}}=0.
\end{array}
\right.
\end{align}
has a unique solution $m$ such that $(\lambda+\tau)^2m \in L^\infty(0,T;L^\infty_\tau(\dot{\mathcal{B}}^{\frac{d}{p}}_{p,1}))$. Moreover, we have
\begin{align}\label{4.7}
\sup_{t\in[0,T]}\|(\lambda+\tau)^2 m\|_{L^\infty_\tau(\dot{\mathcal{B}}^{\frac{d}{p}}_{p,1})}\lm \exp{(\frac{1}{2}\int^T_0F(t')dt')}\|(\lambda+\tau)^2 m_0\|_{L^\infty_\tau(\dot{\mathcal{B}}^{\frac{d}{p}}_{p,1})},
\end{align}
where $F(t)=\frac{1}{\lambda}+1+\|u\|_{\dot{B}^{1+\frac{d}{p}}_{p,1}}+\|u\|^2_{\dot{B}^{1+\frac{d}{p}}_{p,1}}\|S\|^2_{\dot{\mathcal{B}}^{\frac{d}{p}}_{p,1}}$.
\begin{proof}
Multiplying both sides of (\ref{4.6}) by $(\lambda+\tau)^2$, we deduce that
\begin{align}\label{4.8}
\partial_t\widetilde{m}+(u\cdot\nabla)\widetilde{m}+\partial_\tau \widetilde{m}-\frac{2}{\lambda+\tau}\widetilde{m}+(\nabla u:\int^s_0S)\partial_s\widetilde{m}-\partial^2_s\widetilde{m}=0,
\end{align}
where $\widetilde{m}=(\lambda+\tau)^2m$. Applying $\dot{\Delta}_j$ to the both sides of (\ref{4.8}), we obtain that
\begin{align*}
\partial_t\dot{\Delta}_j\widetilde{m}+(u\cdot\nabla)\dot{\Delta}_j\widetilde{m}+\partial_\tau \dot{\Delta}_j\widetilde{m}-\frac{2}{\lambda+\tau}\dot{\Delta}_j\widetilde{m}+\dot{\Delta}_j[(\nabla u:\int^s_0S)\partial_s\widetilde{m}]-\partial^2_s\dot{\Delta}_j\widetilde{m}=R_j,
\end{align*}
where $R_j=[u\nabla, \dot{\Delta}_j]\widetilde{m}$. Multiplying both sides of the above equation by $|\dot{\Delta}_j\widetilde{m}|^{p-1}sign(\dot{\Delta}_j\widetilde{m})$ and integrating over $[-\frac{1}{2},\frac{1}{2}]$ with $s$, we verify that
\begin{align*}
&\frac{1}{p}(\partial_t\int^{\frac{1}{2}}_{-\frac{1}{2}}|\dot{\Delta}_j\widetilde{m}|^{p}ds+(u\cdot\nabla)\int^{\frac{1}{2}}_{-\frac{1}{2}}|\dot{\Delta}_j\widetilde{m}|^{p}ds+\partial_\tau\int^{\frac{1}{2}}_{-\frac{1}{2}}|\dot{\Delta}_j\widetilde{m}|^{p}ds)-\frac{2}{\lambda+\tau}\int^{\frac{1}{2}}_{-\frac{1}{2}}|\dot{\Delta}_j\widetilde{m}|^{p}ds\\
&=\int^{\frac{1}{2}}_{-\frac{1}{2}}\{-\dot{\Delta}_j[(\nabla u:\int^s_0S)\partial_s\widetilde{m}]+R_j+\partial^2_s\dot{\Delta}_j\widetilde{m}\}|\dot{\Delta}_j\widetilde{m}|^{p-1}sign(\dot{\Delta}_j\widetilde{m})ds.
\end{align*}
Integrating by parts and using H\"{o}lder's inequality, we have
\begin{align}\label{4.9}
&\partial_t\|\dot{\Delta}_j\widetilde{m}\|^p_{L^p_s}+(u\cdot\nabla)\|\dot{\Delta}_j\widetilde{m}\|^p_{L^p_s}+\partial_\tau\|\dot{\Delta}_j\widetilde{m}\|^p_{L^p_s}+\int^{\frac{1}{2}}_{-\frac{1}{2}}|\partial_s(\dot{\Delta}_j\widetilde{m})|^2|\dot{\Delta}_j\widetilde{m}|^{p-2}ds
-\frac{1}{\lambda}\|\dot{\Delta}_j\widetilde{m}\|^p_{L^p_s}\\
\nonumber&\lm (\|R_j\|_{L^p_s}+\|\dot{\Delta}_j[(\nabla u:S)\widetilde{m}]\|_{L^p_s})\|\dot{\Delta}_j\widetilde{m}\|^{p-1}_{L^p_s}+\int^{\frac{1}{2}}_{-\frac{1}{2}}\dot{\Delta}_j[(\nabla u:\int^s_0S)\widetilde{m}]\partial_s(\dot{\Delta}_j\widetilde{m})|\dot{\Delta}_j\widetilde{m}|^{p-2}ds.
\end{align}
Using Cauchy-Schwarz's inequality, we get
\begin{multline}\label{4.10}
\int^{\frac{1}{2}}_{-\frac{1}{2}}\dot{\Delta}_j[(\nabla u:\int^s_0S)\widetilde{m}]\partial_s(\dot{\Delta}_j\widetilde{m})|\dot{\Delta}_j\widetilde{m}|^{p-2}ds\\
\leq \varepsilon \int^{\frac{1}{2}}_{-\frac{1}{2}}|\partial_s(\dot{\Delta}_j\widetilde{m})|^2|\dot{\Delta}_j\widetilde{m}|^{p-2}ds+\frac{4}{\varepsilon}\int^{\frac{1}{2}}_{-\frac{1}{2}}|\dot{\Delta}_j[(\nabla u:\int^s_0S)\widetilde{m}]|^2|\dot{\Delta}_j\widetilde{m}|^{p-2}ds,
\end{multline}
where $\varepsilon>0$ is an arbitrary constant. Plugging (\ref{4.10}) into (\ref{4.9}) with $\varepsilon$ small enough, yields that
\begin{align*}
&\partial_t\|\dot{\Delta}_j\widetilde{m}\|^p_{L^p_s}+(u\cdot\nabla)\|\dot{\Delta}_j\widetilde{m}\|^p_{L^p_s}+\partial_\tau\|\dot{\Delta}_j\widetilde{m}\|^p_{L^p_s}+\int^{\frac{1}{2}}_{-\frac{1}{2}}|\partial_s(\dot{\Delta}_j\widetilde{m})|^2|\dot{\Delta}_j\widetilde{m}|^{p-2}ds
-\frac{1}{\lambda}\|\dot{\Delta}_j\widetilde{m}\|^p_{L^p_s}\\
\nonumber&\lm (\|R_j\|_{L^p_s}+\|\dot{\Delta}_j[(\nabla u:S)\widetilde{m}]\|_{L^p_s})\|\dot{\Delta}_j\widetilde{m}\|^{p-1}_{L^p_s}+\int^{\frac{1}{2}}_{-\frac{1}{2}}|\dot{\Delta}_j[(\nabla u:\int^s_0S)\widetilde{m}]|^2|\dot{\Delta}_j\widetilde{m}|^{p-2}ds.
\end{align*}
Integrating the above inequality over $\mathbb{R}^d$ and using the fact that $div~u=0$, we deduce that
\begin{align}\label{4.11}
&\partial_t\|\dot{\Delta}_j\widetilde{m}\|^p_{L^p_{x,s}}+\partial_\tau\|\dot{\Delta}_j\widetilde{m}\|^p_{L^p_{x,s}}+\int_{\mathbb{R}^d}\int^{\frac{1}{2}}_{-\frac{1}{2}}|\partial_s(\dot{\Delta}_j\widetilde{m})|^2|\dot{\Delta}_j\widetilde{m}|^{p-2}dsdx
-\frac{1}{\lambda}\|\dot{\Delta}_j\widetilde{m}\|^p_{L^p_{x,s}}\\
\nonumber&\lm(\|R_j\|_{L^p_{x,s}}+\|\dot{\Delta}_j[(\nabla u:S)\widetilde{m}]\|_{L^p_{x,s}})\|\dot{\Delta}_j\widetilde{m}\|^{p-1}_{L^p_{x,s}}+\|\dot{\Delta}_j[(\nabla u:\int^s_0S)\widetilde{m}]\|^2_{L^p_{x,s}}\|\dot{\Delta}_j\widetilde{m}\|^{p-2}_{L^p_{x,s}},
\end{align}
which implies that
\begin{multline}\label{4.12}
\partial_t\|\dot{\Delta}_j\widetilde{m}\|^2_{L^p_{x,s}}+\partial_\tau\|\dot{\Delta}_j\widetilde{m}\|^2_{L^p_{x,s}}-\frac{1}{\lambda}\|\dot{\Delta}_j\widetilde{m}\|^2_{L^p_{x,s}}\\
\lm (\|R_j\|_{L^p_{x,s}}+\|\dot{\Delta}_j[(\nabla u:S)\widetilde{m}]\|_{L^p_{x,s}})\|\dot{\Delta}_j\widetilde{m}\|_{L^p_{x,s}}+\|\dot{\Delta}_j[(\nabla u:\int^s_0S)\widetilde{m}]\|^2_{L^p_{x,s}}.
\end{multline}
By the definition of Besov spaces and Lemma \ref{L2}, we have
\[\|R_j\|_{L^p_{x,s}}\lm 2^{-j\frac{d}{p}}c_j\|u\|_{\dot{B}^{1+\frac{d}{p}}_{p,1}}\|\widetilde{m}\|_{\dot{\mathcal{B}}^{\frac{d}{p}}_{p,1}}, \]
\[\|\dot{\Delta}_j[(\nabla u:S)\widetilde{m}]\|_{L^p_{x,s}}\lm 2^{-j\frac{d}{p}}c_j\|(\nabla u:S)\widetilde{m}\|_{\dot{\mathcal{B}}^{\frac{d}{p}}_{p,1}}\lm 2^{-j\frac{d}{p}}c_j\|u\|_{\dot{B}^{1+\frac{d}{p}}_{p,1}}\|S\|_{\dot{\mathcal{B}}^{\frac{d}{p}}_{p,1}}\|\widetilde{m}\|_{\dot{\mathcal{B}}^{\frac{d}{p}}_{p,1}},\]
where $c_j\in l^1$. Plugging the above inequality into (\ref{4.12}), we obtain that
\begin{multline}
\partial_t 2^{2j\frac{d}{p}}\|\dot{\Delta}_j\widetilde{m}\|^2_{L^p_{x,s}}+\partial_\tau2^{2j\frac{d}{p}}\|\dot{\Delta}_j\widetilde{m}\|^2_{L^p_{x,s}}
\\
\lm c^2_j\|\widetilde{m}\|^2_{\dot{\mathcal{B}}^{\frac{d}{p}}_{p,1}}\{\frac{1}{\lambda}+\|u\|_{\dot{B}^{1+\frac{d}{p}}_{p,1}}+\|u\|_{\dot{B}^{1+\frac{d}{p}}_{p,1}}\|S\|_{\dot{\mathcal{B}}^{\frac{d}{p}}_{p,1}}+\|u\|^2_{\dot{B}^{1+\frac{d}{p}}_{p,1}}\|S\|^2_{\dot{\mathcal{B}}^{\frac{d}{p}}_{p,1}}\}.
\end{multline}
Taking the $l^1$-norm and using Lemma \ref{L5}, we deduce that
\begin{align}
\sup_{t\in[0,T]}\|\widetilde{m}\|^2_{L^\infty_\tau(\dot{\mathcal{B}}^{\frac{d}{p}}_{p,1})}\lm e^{\displaystyle\int^T_0[\frac{1}{\lambda}+1+\|u\|_{\dot{B}^{1+\frac{d}{p}}_{p,1}}+\|u\|^2_{\dot{B}^{1+\frac{d}{p}}_{p,1}}\|S\|^2_{\dot{\mathcal{B}}^{\frac{d}{p}}_{p,1}}]dt'}(\|\widetilde{m_0}\|^2_{L^\infty_\tau(\dot{\mathcal{B}}^{\frac{d}{p}}_{p,1})}).
\end{align}
From the above estimate and using the theory of the linear transport diffusion equation, one can prove the existence and uniqueness of the equation (\ref{4.6}).
\end{proof}
\end{prop}

\begin{rema}
Under the condition of Proposition \ref{P4.3}. If $\frac{1}{\lambda}$ is small enough, one can obtain a stronger estimate which is independent on $\lambda$.
\end{rema}

\begin{prop}\label{P4.5}
Let $2\leq p<2d$. Assume that $u\in \widetilde{L^2}(0,T;\dot{B}^{1+\frac{d}{p}}_{p,1})$ with $div~u=0$, $w\in \widetilde{L^2}(0,T;\dot{B}^{\frac{d}{p}}_{p,1})$. Suppose that $g,\partial_sg \in L^\infty(0,T;\dot{\mathcal{B}}^{\frac{d}{p}}_{p,1})$, $G,\partial_sG \in L^\infty(0,T;\dot{\mathcal{B}}^{-1+\frac{d}{p}}_{p,1})$, $(\lambda+\tau)^2m_0\in L^\infty_\tau(\dot{\mathcal{B}}^{-1+\frac{d}{p}}_{p,1})$, $(\lambda+\tau)^2M\in L^\infty(0,T; L^\infty_\tau(\dot{\mathcal{B}}^{\frac{d}{p}}_{p,1})$. $m$ is a solution of the following equation 
\begin{align}
\left\{
\begin{array}{ll}
\partial_tm+(u\cdot\nabla)m+\partial_\tau m+g\partial_sm-\partial^2_sm=(w\cdot \nabla) M+G\partial_sM,  \\[1ex]
m|_{t=0}=m_{0}(\tau,x,s), \quad m(t,\tau,x,s)|_{\tau=0}=0, \quad m(t,\tau,x,s)|_{s=\frac{1}{2}}=m(t,\tau,x,s)|_{s=-\frac{1}{2}}=0.
\end{array}
\right.
\end{align}
Let
\[K_0=\sup_{t\in[0,T]}(\|g\|_{\dot{\mathcal{B}}^{\frac{d}{p}}_{p,1}}+\|\partial_sg\|_{\dot{\mathcal{B}}^{\frac{d}{p}}_{p,1}}+\|(\lambda+\tau)^2M\|_{L^\infty_\tau(\dot{\mathcal{B}}^{\frac{d}{p}}_{p,1})}).\]
Then, we have
\begin{align}
\sup_{t\in[0,T]}\|\widetilde{m}\|^2_{L^\infty_\tau(\dot{\mathcal{B}}^{-1+\frac{d}{p}}_{p,1})}\lm e^{\displaystyle\int^T_0(\frac{1}{\lambda}+\|u\|_{\dot{B}^{1+\frac{d}{p}}_{p,1}}+K^2_0)dt'}(\|\widetilde{m_0}\|^2_{L^\infty_\tau(\dot{\mathcal{B}}^{-1+\frac{d}{p}}_{p,1})}+\int^T_0K^2_0(\|w\|^2_{\dot{B}^{\frac{d}{p}}_{p,1}}+\|G\|^2_{\dot{\mathcal{B}}^{-1+\frac{d}{p}}_{p,1}})dt').
\end{align}
\begin{proof}
By the same token as in Proposition \ref{P4.3}, we obtain that
\begin{multline}\label{4.16}
\partial_t\|\dot{\Delta}_j\widetilde{m}\|^2_{L^p_{x,s}}+\partial_\tau\|\dot{\Delta}_j\widetilde{m}\|^2_{L^p_{x,s}}-\frac{1}{\lambda}\|\dot{\Delta}_j\widetilde{m}\|^2_{L^p_{x,s}}-\|\dot{\Delta}_j(G\widetilde{M})\|^2_{L^p_{x,s}}\\
\lm (\|R_j\|_{L^p_{x,s}}+\|\dot{\Delta}_j(w\cdot\nabla \widetilde{M})\|_{L^p_{x,s}}+\|\dot{\Delta}_j(\partial_sg\cdot\widetilde{m})\|_{L^p_{x,s}}+\|\dot{\Delta}_j(\partial_sG\cdot\widetilde{M})\|_{L^p_{x,s}})\|\dot{\Delta}_j\widetilde{m}\|_{L^p_{x,s}}+\|\dot{\Delta}_j(g\widetilde{m})\|^2_{L^p_{x,s}},
\end{multline}
where $\widetilde{m}=(\lambda+\tau)^2m$ and $\widetilde{M}=(\lambda+\tau)^2M$. Note that $2\leq p<2d$.
Taking advantage of Lemmas \ref{L2}-\ref{L3}, we deduce that
\[\|R_j\|_{L^p_{x,s}}\lm 2^{j(1-\frac{d}{p})}c_j\|u\|_{\dot{B}^{1+\frac{d}{p}}_{p,1}}\|\widetilde{m}\|_{\dot{\mathcal{B}}^{-1+\frac{d}{p}}_{p,1}}, \]
\[\|\dot{\Delta}_j(g\cdot\widetilde{m})\|_{L^p_{x,s}}\lm 2^{j(1-\frac{d}{p})}c_j\|g\cdot\widetilde{m}\|_{\dot{\mathcal{B}}^{-1+\frac{d}{p}}_{p,1}}\lm 2^{j(1-\frac{d}{p})}c_j\|g\|_{\dot{\mathcal{B}}^{\frac{d}{p}}_{p,1}}\|\widetilde{m}\|_{\dot{\mathcal{B}}^{-1+\frac{d}{p}}_{p,1}}\lm 2^{j(1-\frac{d}{p})}c_jK_0\|\widetilde{m}\|_{\dot{\mathcal{B}}^{-1+\frac{d}{p}}_{p,1}},\]
\[\|\dot{\Delta}_j(\partial_sg\cdot\widetilde{m})\|_{L^p_{x,s}}\lm 2^{j(1-\frac{d}{p})}c_j\|\partial_sg\cdot\widetilde{m}\|_{\dot{\mathcal{B}}^{-1+\frac{d}{p}}_{p,1}}\lm 2^{j(1-\frac{d}{p})}c_j\|\partial_sg\|_{\dot{\mathcal{B}}^{\frac{d}{p}}_{p,1}}\|\widetilde{m}\|_{\dot{\mathcal{B}}^{-1+\frac{d}{p}}_{p,1}}\lm 2^{j(1-\frac{d}{p})}c_jK_0\|\widetilde{m}\|_{\dot{\mathcal{B}}^{-1+\frac{d}{p}}_{p,1}},\]
\[\|\dot{\Delta}_j(G\cdot\widetilde{M})\|_{L^p_{x,s}}\lm 2^{j(1-\frac{d}{p})}c_j\|G\cdot\widetilde{M}\|_{\dot{\mathcal{B}}^{-1+\frac{d}{p}}_{p,1}}\lm 2^{j(1-\frac{d}{p})}c_j\|G\|_{\dot{\mathcal{B}}^{-1+\frac{d}{p}}_{p,1}}\|\widetilde{M}\|_{\dot{\mathcal{B}}^{\frac{d}{p}}_{p,1}}\lm 2^{j(1-\frac{d}{p})}c_jK_0\|G\|_{\dot{\mathcal{B}}^{-1+\frac{d}{p}}_{p,1}},\]
\[\|\dot{\Delta}_j(\partial_sG\cdot\widetilde{M})\|_{L^p_{x,s}}\lm 2^{j(1-\frac{d}{p})}c_j\|\partial_sG\cdot\widetilde{M}\|_{\dot{\mathcal{B}}^{-1+\frac{d}{p}}_{p,1}}\lm 2^{j(1-\frac{d}{p})}c_j\|\partial_sG\|_{\dot{\mathcal{B}}^{-1+\frac{d}{p}}_{p,1}}\|\widetilde{M}\|_{\dot{\mathcal{B}}^{\frac{d}{p}}_{p,1}}\lm 2^{j(1-\frac{d}{p})}c_jK_0\|G\|_{\dot{\mathcal{B}}^{-1+\frac{d}{p}}_{p,1}},\]
\[\|\dot{\Delta}_j(w\cdot\nabla \widetilde{M})\|_{L^p_{x,s}}\lm 2^{j(1-\frac{d}{p})}c_j\|w\cdot\nabla \widetilde{M}\|_{\dot{\mathcal{B}}^{-1+\frac{d}{p}}_{p,1}}\lm 2^{j(1-\frac{d}{p})}c_j\|w\|_{\dot{B}^{\frac{d}{p}}_{p,1}}\|\widetilde{M}\|_{\dot{\mathcal{B}}^{\frac{d}{p}}_{p,1}}\lm 2^{j(1-\frac{d}{p})}c_jK_0\|w\|_{\dot{B}^{\frac{d}{p}}_{p,1}}, \]
where $c_j\in l^1$. Plugging the above inequality into (\ref{4.16}), we obtain that
\begin{multline}
\partial_t 2^{2j(-1+\frac{d}{p})}\|\dot{\Delta}_j\widetilde{m}\|^2_{L^p_{x,s}}+\partial_\tau2^{2j(-1+\frac{d}{p})}\|\dot{\Delta}_j\widetilde{m}\|^2_{L^p_{x,s}}
\\
\lm c^2_j\{\|\widetilde{m}\|^2_{\dot{\mathcal{B}}^{\frac{d}{p}}_{p,1}}[\frac{1}{\lambda}+\|u\|_{\dot{B}^{1+\frac{d}{p}}_{p,1}}+K^2_0]+K^2_0(\|w\|^2_{\dot{B}^{\frac{d}{p}}_{p,1}}+\|G\|^2_{\dot{\mathcal{B}}^{-1+\frac{d}{p}}_{p,1}})\}.
\end{multline}
Taking the $l^1$-norm and using Lemma \ref{L5}, we deduce that
\begin{align}
\sup_{t\in[0,T]}\|\widetilde{m}\|^2_{L^\infty_\tau(\dot{\mathcal{B}}^{-1+\frac{d}{p}}_{p,1})}\lm e^{\displaystyle\int^T_0(\frac{1}{\lambda}+\|u\|_{\dot{B}^{1+\frac{d}{p}}_{p,1}}+K^2_0)dt'}(\|\widetilde{m_0}\|^2_{L^\infty_\tau(\dot{\mathcal{B}}^{-1+\frac{d}{p}}_{p,1})}+\int^T_0K^2_0(\|w\|^2_{\dot{B}^{\frac{d}{p}}_{p,1}}+\|G\|^2_{\dot{\mathcal{B}}^{-1+\frac{d}{p}}_{p,1}})dt').
\end{align}
\end{proof}
\end{prop}

\subsection{A priori estimates for the linear Navier-Stokes equations}
In this subsection, we review a priori estimates for the linear Navier-Stokes equations.
\begin{lemm}\cite{Danchin2007}\label{L4.6}
Let $1\leq p\leq \infty$, $1\leq \rho_1\leq \rho\leq \infty$ and $-1-d\min\{\frac{1}{p},\frac{1}{p'}\}<s\leq 1+\frac{d}{p}$ where $p'$ is the conjugate number of $p$. Assume that $v_{0}\in \dot{B}^{s}_{p,1}$, $u\in \widetilde{L}^2(0,T;\dot{B}^{1+\frac{d}{p}}_{p,1})\hookrightarrow L^{1}(0,T;\dot{B}^{1+\frac{d}{p}}_{p,1})$ with $div~u=0$, $f\in\widetilde{L}^2_T(\dot{B}^{s-1}_{p,1})$. Suppose that $v$ is the solution of
\begin{align}\label{4.20}
\left\{
\begin{array}{ll}
\partial_{t}v+(u\cdot\nabla)v-(1-\omega)\Delta v+\nabla{P}=f, ~ div~v=0,\\[1ex]
v|_{t=0}=v_{0} .\\[1ex]
\end{array}
\right.
\end{align}
Then, we have following estimates:
$$\|v\|_{\widetilde{L}^{\rho}_T(\dot{B}^{s+\frac{2}{\rho}}_{p,1})}\leq Ce^{CU(T)}(\|v_{0}\|_{\dot{B}^{s}_{p,1}}+\|f\|_{\widetilde{L}^{\rho_1}_T(\dot{B}^{s-2+\frac{2}{\rho_1}}_{p,1})}),$$
$$\|\nabla P\|_{\widetilde{L}^{\rho}_T(\dot{B}^{s+\frac{2}{\rho}}_{p,1})}\leq (e^{CU(T)}-1)(\|v_{0}\|_{\dot{B}^{s}_{p,1}}+\|f\|_{\widetilde{L}^{\rho_1}_T(\dot{B}^{s-2+\frac{2}{\rho_1}}_{p,1})}),$$
where $U(t)=\int^{t}_{0}\|\nabla u\|_{\dot{B}^{\frac{d}{p}}_{p,1}}dt'$.
\end{lemm}

\begin{coro}
Let $2\leq p<2d$ where $p'$. Assume that $v_{0}\in \dot{B}^{\frac{d}{p}}_{p,1}$, $u\in \widetilde{L}^\infty(0,T;\dot{B}^{\frac{d}{p}}_{p,1})\cap \widetilde{L}^2(0,T;\dot{B}^{1+\frac{d}{p}}_{p,1}) $ with $div~u=div~v_0=0$, $f\in\widetilde{L}^2_T(\dot{B}^{\frac{d}{p}-1}_{p,1})$. Suppose that $v$ is the solutions of (\ref{4.20}) and $v_L$ is the solution of
\begin{align}\label{4.21}
\left\{
\begin{array}{ll}
\partial_{t}v_L-(1-\omega)\Delta v_L=0,\\[1ex]
v_L|_{t=0}=v_{0}. \\[1ex]
\end{array}
\right.
\end{align}
Then, we have the following estimate:
$$\|v\|_{\widetilde{L}^2_T(\dot{B}^{\frac{d}{p}+1}_{p,1})}\leq \|e^{t(1-\omega)\Delta}v_0\|_{\widetilde{L}^2_T(\dot{B}^{\frac{d}{p}+1}_{p,1})}+Ce^{CU(T)}(\|f\|_{\widetilde{L}^2_T(\dot{B}^{\frac{d}{p}-1}_{p,1})}+T^{\frac{1}{2}}\|u\|_{\widetilde{L}^\infty_T(\dot{B}^{\frac{d}{p}}_{p,1})}\| v_0\|_{\dot{B}^{\frac{d}{p}}_{p,1}}).$$
\begin{proof}
Let $\overline{v}\doteq v-v_L$. By the equations (\ref{4.20}) and (\ref{4.21}), we see that $\overline{v}$ satisfies that
\begin{align}\label{4.22}
\left\{
\begin{array}{ll}
\partial_{t}\overline{v}+(u\cdot\nabla)\overline{v}-(1-\omega)\Delta \overline{v}+\nabla{P}=f-(u\cdot \nabla) v_L, ~ div~\overline{v}=0,\\[1ex]
\overline{v}|_{t=0}=0 .\\[1ex]
\end{array}
\right.
\end{align}
Taking advantage of Lemma \ref{L4.6} with $s=\frac{d}{p}$ and $\rho=\rho_1=2$, we deduce that
\begin{align}\label{4.23}
\|\overline{v}\|_{\widetilde{L}^2_T(\dot{B}^{\frac{d}{p}+1}_{p,1})}\leq Ce^{CU(T)}(\|f\|_{\widetilde{L}^2_T(\dot{B}^{\frac{d}{p}-1}_{p,1})}+\|(u\cdot\nabla) v_L\|_{\widetilde{L}^2_T(\dot{B}^{\frac{d}{p}-1}_{p,1})}).
\end{align}
Using Lemma \ref{L3} with $s_1=\frac{d}{p}$ and $s_2=\frac{d}{p}-1$, we obtain that
$$\|(u\cdot\nabla) v_L\|_{\widetilde{L}^2_T(\dot{B}^{\frac{d}{p}-1}_{p,1})}\lm T^{\frac{1}{2}}\|u\|_{\widetilde{L}^\infty_T(\dot{B}^{\frac{d}{p}}_{p,1})}\|\nabla v_L\|_{\widetilde{L}^\infty_T(\dot{B}^{\frac{d}{p}-1}_{p,1})}\lm T^{\frac{1}{2}}\|u\|_{\widetilde{L}^\infty_T(\dot{B}^{\frac{d}{p}}_{p,1})}\| v_L\|_{\widetilde{L}^\infty_T(\dot{B}^{\frac{d}{p}}_{p,1})}.$$
\end{proof}
Plugging the above inequality into (\ref{4.23}), yields that
\begin{align}
\|\overline{v}\|_{\widetilde{L}^2_T(\dot{B}^{\frac{d}{p}+1}_{p,1})}\leq Ce^{CU(T)}(\|f\|_{\widetilde{L}^2_T(\dot{B}^{\frac{d}{p}-1}_{p,1})}+T^{\frac{1}{2}}\|u\|_{\widetilde{L}^\infty_T(\dot{B}^{\frac{d}{p}}_{p,1})}\| v_L\|_{\widetilde{L}^\infty_T(\dot{B}^{\frac{d}{p}}_{p,1})}),
\end{align}
which leads to
$$\|v\|_{\widetilde{L}^2_T(\dot{B}^{\frac{d}{p}+1}_{p,1})}\leq \|v_L\|_{\widetilde{L}^2_T(\dot{B}^{\frac{d}{p}+1}_{p,1})}+Ce^{CU(T)}(\|f\|_{\widetilde{L}^2_T(\dot{B}^{\frac{d}{p}-1}_{p,1})}+T^{\frac{1}{2}}\|u\|_{\widetilde{L}^\infty_T(\dot{B}^{\frac{d}{p}}_{p,1})}\| v_L\|_{\widetilde{L}^\infty_T(\dot{B}^{\frac{d}{p}}_{p,1})}).$$
Using the fact that $v_L=e^{t(1-\omega)\Delta}v_0$ and $\|e^{t(1-\omega)\Delta}v_0\|_{\widetilde{L}^\infty_T(\dot{B}^{\frac{d}{p}}_{p,1})}\lm \|v_0\|_{\dot{B}^{\frac{d}{p}}_{p,1}}$ then completes the proof.
\end{coro}

\section{Local well-posedness}
\subsection{Approximate solutions}
~~~~~First, we construct approximate solutions which are smooth (for $x$ variable) solutions of some linear equations.\\ [1ex]
~~~~~Starting for $(v^{0},F^{0}, m^0)\triangleq (\dot{S}_{0}v_{0},\dot{S}_{0}F_{0}, \dot{S}_0m_0)$,  we define by induction a sequence $(v^{n},F^{n}, m^n, P^n)_{n\in\mathbb{N}}$  by solving the following linear equations:
\begin{align}\tag{5.1n}\label{5.1n}
\left\{
\begin{array}{ll}
\partial_{t}v^{n+1}+(v^{n}\cdot\nabla)v^{n+1}-(1-\omega)\Delta v^{n+1}+\nabla{P^{n+1}}=div~\sigma^{n},  \quad div~ v^{n+1}=0,\\[1ex]
\sigma^{n}=\omega\int^{\frac{1}{2}}_{-\frac{1}{2}}S^nds,\quad S^n=\displaystyle\int^{+\infty}_0m^n\widetilde{\mathscr{S}}(F^n)d\tau, \\[1ex]
\partial_{t}F^{n+1}+\partial_\tau F^{n+1}+(v^{n}\cdot\nabla)F^{n+1}=F^{n+1}\cdot\nabla v^n+\nabla v^n,  \\[1ex]
\partial_tm^{n+1}+(v^n\cdot\nabla)m^{n+1}+\partial_\tau m^{n+1}+(\nabla v^n:\int^s_0S^n)\partial_sm^{n+1}-\partial^2_sm^{n+1}=0, \\[1ex]
v^{n+1}|_{t=0}=\dot{S}_{n+1}v_{0},\quad F^{n+1}|_{t=0}=\dot{S}_{n+1}F_{0} ,\quad m^{n+1}|_{t=0}=\dot{S}_{n+1}m_0, \\[1ex]
m^{n+1}|_{\tau=0}=0,\quad m^{n+1}|_{s=\frac{1}{2}}=m^{n+1}|_{s=-\frac{1}{2}}=0, \quad  F^{n+1}|_{\tau=0}=0.\\[1ex]
\end{array}
\right.
\end{align}
For a fixed $n\geq 1$, assume that $(v^{n},F^{n}, m^n)$ is smooth enough. Since the initial data are smooth, it follows that there exists a smooth solution $(v^{n+1},F^{n+1}, m^{n+1},P^{n+1})$ satisfies the linear equations (\ref{5.1n}).

\subsection{Uniform bounds}
Next, we are going to find some positive $T$ such that for which the approximate solutions are uniformly bounded. Define $V^{n}(t)=\displaystyle\int^{t}_{0}\|v^n(t')\|_{\dot{B}^{1+\frac{d}{p}}_{p,1}}dt'$. Then by Propositions \ref{P4.1}, \ref{P4.3} and Lemma \ref{4.6} with $f=div~\sigma^n$, we have
\begin{align}\label{5.1}
\|F^{n+1}\|_{L^\infty_\tau(\dot{B}^{\frac{d}{p}}_{p,1})}\leq C \exp(V^n(t))\bigg[\|F_0\|_{L^\infty_\tau(\dot{B}^{\frac{d}{p}}_{p,1})}+\int^t_0\|v^n(t')\|_{\dot{B}^{1+\frac{d}{p}}_{p,1}}dt'\bigg],
\end{align}
\begin{align*}
\|(\lambda+\tau)^2 m^{n+1}\|_{L^\infty_\tau(\dot{\mathcal{B}}^{\frac{d}{p}}_{p,1})}\leq C\exp{\int^t_0\frac{1}{2}(\frac{1+\lambda}{\lambda}+\|v^n\|_{\dot{B}^{1+\frac{d}{p}}_{p,1}}+\|v^n\|^2_{\dot{B}^{1+\frac{d}{p}}_{p,1}}\|S^n\|^2_{\dot{\mathcal{B}}^{\frac{d}{p}}_{p,1}})dt'}\|(\lambda+\tau)^2 m_0\|_{L^\infty_\tau(\dot{\mathcal{B}}^{\frac{d}{p}}_{p,1})},
\end{align*}
\begin{align*}
\|v^{n+1}\|_{\widetilde{L}^\infty_t(B^{\frac{d}{p}}_{p,1})}\leq C\exp{CV^n(t)}(\|v_{0}\|_{\dot{B}^{\frac{d}{p}}_{p,1}}+\|\sigma^n\|_{\widetilde{L}^2_t(\dot{B}^{\frac{d}{p}}_{p,1})}),
\end{align*}
\begin{align*}
\|v^{n+1}\|_{\widetilde{L}^2_t(\dot{B}^{\frac{d}{p}+1}_{p,1})}\leq \|e^{t'(1-\omega)\Delta}v_0\|_{\widetilde{L}^2_t(\dot{B}^{\frac{d}{p}+1}_{p,1})}+Ce^{CV^n(t)}(\|\sigma^n\|_{\widetilde{L}^2_t(\dot{B}^{\frac{d}{p}}_{p,1})}+t^{\frac{1}{2}}\|v^n\|_{\widetilde{L}^\infty_t(\dot{B}^{\frac{d}{p}}_{p,1})}\| v_0\|_{\dot{B}^{\frac{d}{p}}_{p,1}})
\end{align*}
And by Corollary \ref{C2}, we obtain
\begin{align}\label{5.2}
\|\sigma^n\|_{\widetilde{L}^2_t(\dot{B}^{\frac{d}{p}}_{p,1})}\leq Ct^{\frac{1}{2}}\|S^n\|_{\dot{\mathcal{B}}^{\frac{d}{p}}_{p,1}}\leq\frac{Ct^{\frac{1}{2}}}{\lambda}\|(\lambda+\tau)^2m^{n}\|_{L^\infty_\tau(\dot{B}^{\frac{d}{p}}_{p,1})}\|F^n\|_{L^\infty_\tau(\dot{B}^{\frac{d}{p}}_{p,1})} .
\end{align}
Let $R_0\doteq\|F_0\|_{L^\infty_\tau(\dot{B}^{\frac{d}{p}}_{p,1})}+\|v_{0}\|_{\dot{B}^{\frac{d}{p}}_{p,1}}+\|(\lambda+\tau)^2 m_0\|^2_{L^\infty_\tau(\dot{\mathcal{B}}^{\frac{d}{p}}_{p,1})}$.
We define that
\[R(t)=\frac{4CR_0}{1-CtR_0}, \quad  A(t)=\|e^{t'(1-\omega)\Delta}v_0\|_{\widetilde{L}^2_t(\dot{B}^{\frac{d}{p}+1}_{p,1})}+2Ct^{\frac{1}{2}}(\frac{R^2(t)}{\lambda}+R(t)R_0).\]

Note that $R(t)\rightarrow 4CR_0$. Since $\|e^{t'(1-\omega)\Delta}v_0\|_{\widetilde{L}^2_t(\dot{B}^{\frac{d}{p}+1}_{p,1})}\rightarrow 0 $ as $t\rightarrow 0$(See Chapter 10 in \cite{Bahouri2011}), it follows that
$A(t)\rightarrow 0$ as $t\rightarrow 0$, which implies that
$$ e^{[t(\frac{1+\lambda}{\lambda})+t^{\frac{1}{2}}A(t)+CA^2(t)R^4(t)]} \rightarrow 1, \quad  t\rightarrow 0.$$
Thus, we can find a $T_1>0$ such that
\[e^{[t(\frac{1+\lambda}{\lambda})+t^{\frac{1}{2}}A(t)+CA^2(t)R^4(t)]}\leq 2, \quad \forall t\in[0,T_1]. \]
On the other hand, for any fixed $\lambda$, there exists a $T_2>0$, such that \[\frac{t^{\frac{1}{2}}R(t)}{\lambda}<R_0,  \quad t^{\frac{1}{2}}A(t)\leq R_0.\]

Taking $T=\min\{T_1, T_2, \frac{1}{2CR_0}\}$, we obtain that for any $t\in [0, T]$
\[e^{[t(\frac{1+\lambda}{\lambda})+t^{\frac{1}{2}}A(T)+CA^2(t)R^4(t)]}\leq 2, \quad 4CR_0\leq R(t)\leq 8CR_0, \]
\[\frac{t^{\frac{1}{2}}R(t)}{\lambda}\leq R_0, \quad t^{\frac{1}{2}}A(t)\leq R_0.\]

Suppose that for any $t\in[0,T]$
$$\|F^{n}\|_{L^\infty_\tau(\dot{B}^{\frac{d}{p}}_{p,1})}\leq R(t), \quad \|(\lambda+\tau)^2 m^{n}\|^2_{L^\infty_\tau(\dot{\mathcal{B}}^{\frac{d}{p}}_{p,1})}\leq R(t),\quad  \|v^n\|_{\widetilde{L}^\infty_t(B^{\frac{d}{p}}_{p,1})}\leq R(t),
\quad \|v^n\|_{\widetilde{L}^2_t(\dot{B}^{\frac{d}{p}+1}_{p,1})}\leq A(t).$$

By virtue of (\ref{5.1}) and (\ref{5.2}), we see that
\begin{align}\label{5.3}
\sup_{t\in[0,T]}\|F^{n+1}\|_{L^\infty_\tau(\dot{B}^{\frac{d}{p}}_{p,1})}\leq Ce^{T^{\frac{1}{2}}A(T)}(R_0+T^{\frac{1}{2}}A(T))\leq 2C(R_0+T^{\frac{1}{2}}R(T))\leq 4CR_0\leq R(t), \\
\sup_{t\in[0,T]}\|(\lambda+\tau)^2 m^{n+1}\|_{L^\infty_\tau(\dot{\mathcal{B}}^{\frac{d}{p}}_{p,1})}\leq Ce^{\frac{1}{2}[T(\frac{1+\lambda}{\lambda})+T^{\frac{1}{2}}A(T)+CA^2(T)R^4(T)]}R_0\leq 2CR_0\leq R(t),\\\label{5.5}
\|v^{n+1}\|_{\widetilde{L}^\infty_T(B^{\frac{d}{p}}_{p,1})}\leq Ce^{T^{\frac{1}{2}}A(T)}(R_0+\frac{T^{\frac{1}{2}}}{\lambda}R^2(T))\leq 4CR_0\leq R(t),\\
\|v^{n+1}\|_{\widetilde{L}^2_t(B^{\frac{d}{p}+1}_{p,1})}\leq \|e^{-t'(1-\omega)\Delta}v_0\|_{\widetilde{L}^2_t(\dot{B}^{\frac{d}{p}+1}_{p,1})}+Ce^{t^{\frac{1}{2}}A(t)}t^{\frac{1}{2}}(R_0R(t)+\frac{ R^2(t)}{\lambda})\leq A(t).
\end{align}

 By induction, we obtain that $\|v^n\|_{\widetilde{L}^\infty_T(B^{\frac{d}{p}+1}_{p,1})\cap \widetilde{L}^2_T(B^{\frac{d}{p}+1}_{p,1})}$, $\sup_{t\in[0,T]}\|F^n\|_{L^\infty_\tau(\dot{B}^{\frac{d}{p}}_{p,1})}$ and $\sup_{t\in[0,T]}\|(\lambda+\tau)^2 m^{n}\|_{L^\infty_\tau(\dot{\mathcal{B}}^{\frac{d}{p}}_{p,1})}$ is uniformly bounded.  Taking advantage of Lemma \ref{L4.6}, we deduce that $\|\nabla P^n\|_{\widetilde{L}^\infty_T(B^{\frac{d}{p}+1}_{p,1})}$ is also uniformly bounded.

\subsection{Convergence}
 We are going to show that $(v^{n},F^{n}, m^n)_{n\in \mathbb{N}}$ is a Cauchy sequence in some Banach spaces. Let $\overline{v}^{n,l}=v^{n}-v^{l}, \overline{F}^{n,l}=F^{n}-F^{l},\overline{m}^{n,l}=m^{n}-m^{l}$ with $1\leq l<n$. By (\ref{5.1n}) we have
 \begin{align}\label{5.7}
\left\{
\begin{array}{ll}
\partial_{t}\overline{v}^{n,l}+(v^{n-1}\cdot\nabla)(\overline{v}^{n,l})-(1-\omega)\Delta \overline{v}^{n,l}+\nabla(P^{n}-P^{l})=div~(\sigma^{n-1}-\sigma^{l-1})-\overline{v}^{n-1,l-1}\nabla v^l,  \\[1ex]
\partial_{t}\overline{F}^{n,l}+(v^{n-1}\cdot\nabla)\overline{F}^{n,l}+\partial_\tau \overline{F}^{n,l}=-\overline{F}^{n-1,l-1}\cdot \nabla v^l+\overline{F}^{n,l}\cdot \nabla v^{l-1}+F^l\nabla\overline{v}^{n-1,l-1}+\nabla\overline{v}^{n-1,l-1},  \\[1ex]
\partial_{t}\overline{m}^{n,l}+(v^{n-1}\cdot\nabla)\overline{m}^{n,l}+\partial_\tau \overline{m}^{n,l}-\partial^2_s\overline{m}^{n,l}+(\nabla v^{n-1}:\int^s_0S^{n-1})\partial_s\overline{m}^{n,l}=-\overline{v}^{n-1,l-1}\nabla m^l+H^{n,l}\partial_sm^l,  \\[1ex]
\overline{v}^{n,l}|_{t=0}=(S_n-S_l)v_{0}, \quad \overline{F}^{n,l}|_{t=0}=(S_n-S_l)F_{0},\quad \overline{m}^{n,l}|_{t=0}=(S_n-S_l)m_{0}, \\[1ex]
\overline{F}^{n,l}|_{\tau=0}=0,\quad \overline{m}^{n,l}|_{\tau=0}=0, \quad \overline{m}^{n,l}|_{s=\pm\frac{1}{2}}=0,
\end{array}
\right.
\end{align}
where $H^{n,l}=\nabla v^{l-1}:\int^s_0S^{l-1}-\nabla v^{n-1}:\int^s_0S^{n-1}$. Note that
\[\|F^{n}\|_{L^\infty_\tau(\dot{B}^{\frac{d}{p}}_{p,1})},~
\|(\lambda+\tau)^2 m^{n}\|_{L^\infty_\tau(\dot{\mathcal{B}}^{\frac{d}{p}}_{p,1})},~
\|v^{n}\|_{\widetilde{L}^\infty_T(B^{\frac{d}{p}}_{p,1})\cap \widetilde{L}^2_T(\dot{B}^{\frac{d}{p}+1}_{p,1})}\leq R(t)\leq 8CR_0\].
By virtue of Propositions \ref{P4.2}, \ref{P4.5} and Lemma \ref{L4.6}, we deuce that
\begin{align}\label{5.7}
\nonumber\sup_{t\in[0,T]}\|\overline{F}^{n,l}\|_{L^\infty_\tau(\dot{B}^{-1+\frac{d}{p}}_{p,1})}&\lm \|(S_n-S_l)F_0\|_{L^\infty_\tau(\dot{B}^{-1+\frac{d}{p}}_{p,1})}+\int^T_0\|\overline{v}^{n-1,l-1}\|_{\dot{B}^{\frac{d}{p}}_{p,1}}(1+\|F^l\|_{\dot{B}^{\frac{d}{p}}_{p,1}})dt'\\
&\lm A^1_{n,l}+\int^T_0\|\overline{v}^{n-1,l-1}\|_{\dot{B}^{\frac{d}{p}}_{p,1}}dt',
\end{align}
\begin{align}\label{5.8}
\nonumber\sup_{t\in[0,T]}\|(\lambda+\tau)^2\overline{m}^{n,l}\|^2_{L^\infty_\tau(\dot{\mathcal{B}}^{-1+\frac{d}{p}}_{p,1})}&\lm\|(\lambda+\tau)^2(S_n-S_l){m_0}\|^2_{L^\infty_\tau(\dot{\mathcal{B}}^{-1+\frac{d}{p}}_{p,1})}+\int^T_0(\|\overline{v}^{n-1,l-1}\|^2_{\dot{B}^{\frac{d}{p}}_{p,1}}+\|G^{n,l}\|^2_{\dot{\mathcal{B}}^{-1+\frac{d}{p}}_{p,1}})dt'\\
&\lm A^2_{n,l}+\int^T_0(\|\overline{v}^{n-1,l-1}\|^2_{\dot{B}^{\frac{d}{p}}_{p,1}}+\|G^{n,l}\|^2_{\dot{\mathcal{B}}^{-1+\frac{d}{p}}_{p,1}})dt'),
\end{align}
\begin{align}\label{5.9}
\nonumber\|\overline{v}^{n,l}\|_{\widetilde{L}^\infty_T(\dot{B}^{-1+\frac{d}{p}}_{p,1})\cap \widetilde{L}^2_T(\dot{B}^{\frac{d}{p}}_{p,1})}&\lm(\|(S_n-S_l)v_{0}\|_{\dot{B}^{-1+\frac{d}{p}}_{p,1}}+\|div~(\sigma^{n-1}-\sigma^{l-1})\|_{\widetilde{L}^2_T(\dot{B}^{-2+\frac{d}{p}}_{p,1})}+\|\overline{v}^{n-1,l-1}\nabla v^l\|_{\widetilde{L}^2_T(\dot{B}^{-2+\frac{d}{p}}_{p,1})})\\
&\lm A^3_{n,l}+\|div~(\sigma^{n-1}-\sigma^{l-1})\|_{\widetilde{L}^2_T(\dot{B}^{-2+\frac{d}{p}}_{p,1})}+\|\overline{v}^{n-1,l-1}\nabla v^l\|_{\widetilde{L}^2_T(\dot{B}^{-2+\frac{d}{p}}_{p,1})}),
\end{align}
where $A^1_{n,l}, A^2_{n,l}, A^3_{n,l}\rightarrow 0$ as $n,l\rightarrow \infty$(Here, we use the fact that $v_0\in \dot{B}^{-1+\frac{d}{p}}_{p,1}$, $F_0\in L^\infty_\tau(\dot{B}^{-1+\frac{d}{p}}_{p,1})$ and $(\lambda+\tau)^2m_0\in L^\infty_\tau(\dot{B}^{-1+\frac{d}{p}}_{p,1})$).

Taking advantage of Lemma \ref{L3}, we obtain
\begin{align}\label{5.10}
&\|div~(\sigma^{n-1}-\sigma^{l-1})\|_{\widetilde{L}^2_T(\dot{B}^{-2+\frac{d}{p}}_{p,1})}\lm \|\sigma^{n-1}-\sigma^{l-1}\|_{\widetilde{L}^2_T(\dot{B}^{-1+\frac{d}{p}}_{p,1})}\lm T^{\frac{1}{2}}\|S^{n-1}-S^{l-1}\|_{L^\infty_T(\dot{\mathcal{B}}^{-1+\frac{d}{p}}_{p,1})}\\
\nonumber &\lm \frac{T^{\frac{1}{2}}}{\lambda}[\|(\lambda+\tau)^2\overline{m}^{n-1,l-1}\|_{L^\infty_\tau(\dot{\mathcal{B}}^{-1+\frac{d}{p}}_{p,1})}\|F^{l-1}\|_{L^\infty_\tau(\dot{B}^{\frac{d}{p}}_{p,1})}+\|(\lambda+\tau)^2\overline{m}^{n-1}\|_{L^\infty_\tau(\dot{\mathcal{B}}^{\frac{d}{p}}_{p,1})}\|\overline{F}^{n-1,l-1}\|_{L^\infty_\tau(\dot{B}^{-1+\frac{d}{p}}_{p,1})}]\\
\nonumber &\lm\frac{T^{\frac{1}{2}}}{\lambda}[\|(\lambda+\tau)^2\overline{m}^{n-1,l-1}\|_{L^\infty_\tau(\dot{\mathcal{B}}^{-1+\frac{d}{p}}_{p,1})}+\|\overline{F}^{n-1,l-1}\|_{L^\infty_\tau(\dot{B}^{-1+\frac{d}{p}}_{p,1})}].
\end{align}
Since $div~v^n=0$, it follows that
\begin{multline}\label{5.11}
\|\overline{v}^{n-1,l-1}\nabla v^l\|_{\widetilde{L}^2_T(\dot{B}^{-2+\frac{d}{p}}_{p,1})}=\|div~(\overline{v}^{n-1,l-1}\otimes v^l)\|_{\widetilde{L}^2_T(\dot{B}^{-2+\frac{d}{p}}_{p,1})}\\
\lm\|\overline{v}^{n-1,l-1}\otimes v^l\|_{\widetilde{L}^2_T(\dot{B}^{-1+\frac{d}{p}}_{p,1})}\lm \|\overline{v}^{n-1,l-1}\|_{\widetilde{L}^\infty_T(\dot{B}^{-1+\frac{d}{p}}_{p,1})}\|v^l\|_{\widetilde{L}^2_T(\dot{B}^{\frac{d}{p}}_{p,1})}\lm V(T)\|\overline{v}^{n-1,l-1}\|_{\widetilde{L}^\infty_T(\dot{B}^{-1+\frac{d}{p}}_{p,1})},
\end{multline}
where $V(T)=\sup_{l}\|v^l\|_{\widetilde{L}^2_T(\dot{B}^{\frac{d}{p}}_{p,1})}$. Notice that $H^{n,l}=\nabla v^{l-1}:\int^s_0S^{l-1}-\nabla v^{n-1}:\int^s_0S^{n-1}$. A standard triangle inequality leads to
\begin{multline}\label{5.12}
\|H^{n,l}\|_{L^2_T\dot{\mathcal{B}}^{-1+\frac{d}{p}}_{p,1}}\lm V(T)[\|(\lambda+\tau)^2\overline{m}^{n-1,l-1}\|_{L^\infty_\tau(\dot{\mathcal{B}}^{-1+\frac{d}{p}}_{p,1})}+\|\overline{F}^{n-1,l-1}\|_{L^\infty_\tau(\dot{B}^{-1+\frac{d}{p}}_{p,1})}]+\|\overline{v}^{n-1,l-1}\|_{\widetilde{L}^2_T(\dot{B}^{\frac{d}{p}}_{p,1})}.
\end{multline}

Assume that $$X(T)=\limsup_{n,l\rightarrow\infty}\sup_{t\in[0,T]}\|\overline{F}^{n,l}\|_{L^\infty_\tau(\dot{B}^{-1+\frac{d}{p}}_{p,1})},\quad Y(T)=\limsup_{n,l\rightarrow\infty}\sup_{t\in[0,T]}\|(\lambda+\tau)^2\overline{m}^{n,l}\|^2_{L^\infty_\tau(\dot{\mathcal{B}}^{-1+\frac{d}{p}}_{p,1})},$$ $$Z(T)=\limsup_{n,l\rightarrow\infty}\|\overline{v}^{n,l}\|_{\widetilde{L}^\infty_T(\dot{B}^{-1+\frac{d}{p}}_{p,1})\cap \widetilde{L}^2_T(\dot{B}^{\frac{d}{p}}_{p,1})}.$$
Plugging (\ref{5.10})-(\ref{5.12}) into (\ref{5.7})-(\ref{5.9}) and taking the limit as $n,l$ go to infinite, we deduce that
\[ X(T)\lm T^{\frac{1}{2}} Z(T), \quad Y(T)\lm V(T)(X(T)+(Y(T)))+Z(T), \]
\[Z(T)\lm T^{\frac{1}{2}}(X(T)+Y(T))+V(T)Z(T),\]
which leads to
\[X(T)+\varepsilon Y(T)+Z(T)\lm(V(T)+T^{\frac{1}{2}})X(T)+(\varepsilon V(T)+T^{\frac{1}{2}})Y(T)+(\varepsilon+T^{\frac{1}{2}})Z(T).\]
Since $V(T)\rightarrow 0$ as $ T\rightarrow 0$, it follows that we can choose $\varepsilon$ and $T$ small enough such that
\[X(T)+\varepsilon Y(T)+Z(T)\leq \frac{1}{2}(X(T)+\varepsilon Y(T)+Z(T)).\]
Hence, $X(T), Y(T), Z(T)=0$ and then we deduce that $(v^{n},F^{n}, m^n)_{n\in \mathbb{N}}$ is a Cauchy sequence. Moreover, there exists a limit function $(v, F, m)$ such that $v^n\rightarrow v$ in $\widetilde{L^\infty}(B^{-1+\frac{d}{p}}_{p,1})$,  $F^n\rightarrow F$ in $L^\infty_{t,\tau}(\dot{B}^{-1+\frac{d}{p}}_{p,1})$ and $(\lambda+\tau)^2m^n\rightarrow (\lambda+\tau)^2m$ in $L^\infty_{t,\tau}(\dot{\mathcal{B}}^{-1+\frac{d}{p}}_{p,1})$.
By virtue of the equation (\ref{5.1n}), one can obtain $P^n$ by solving an elliptic equation. Indeed, $v^n$ is a Cauchy sequence implies that $P_n$ is also a Cauchy sequence and there exists a limit function $P$ such that $P^n \rightarrow P$. (See, \cite{Bahouri2011,Danchin2007})

Passing to the limit in (\ref{5.1n}) in the weak sense, we conclude that  $(v,F,m,P)$ is indeed a solution of (\ref{1.5}). Since
\[\sup_{t\in[0,T]}\|F^{n}\|_{L^\infty_\tau(\dot{B}^{\frac{d}{p}}_{p,1})},
\sup_{t\in[0,T]}\|(\lambda+\tau)^2 m^{n}\|^2_{L^\infty_\tau(\dot{\mathcal{B}}^{\frac{d}{p}}_{p,1})},
\|v^{n}\|_{\widetilde{L}^\infty_T(B^{\frac{d}{p}}_{p,1})\cap \widetilde{L}^2_T(\dot{B}^{\frac{d}{p}+1}_{p,1})}\leq 8CR_0,\]
it follow by the Fatou properties of Besov spaces that
\[\sup_{t\in[0,T]}\|F\|_{L^\infty_\tau(\dot{B}^{\frac{d}{p}}_{p,1})},
\sup_{t\in[0,T]}\|(\lambda+\tau)^2 m\|^2_{L^\infty_\tau(\dot{\mathcal{B}}^{\frac{d}{p}}_{p,1})},
\|v\|_{\widetilde{L}^\infty_T(B^{\frac{d}{p}}_{p,1})\cap \widetilde{L}^2_T(\dot{B}^{\frac{d}{p}+1}_{p,1})}\leq 8CR_0.\]
Using the equations, one can see that $\partial_t v, \partial_t F, (\lambda+\tau)^2\partial_tm$ belongs to some Besov spaces with negative regularity, which implies that
$$v\in C([0,T^{*});\dot{B}^{-1+\frac{d}{p}}_{p,1}\cap \dot{B}^{\frac{d}{p}}_{p,1}, $$
$$ F\in C([0,T^{*});L^\infty_\tau(\dot{B}^{-1+\frac{d}{p}}_{p,1})\cap \dot{B}^{\frac{d}{p}}_{p,1})),~~~(\lambda+\tau)^2m\in C([0,T^{*});L^\infty_\tau(\dot{\mathcal{B}}^{-1+\frac{d}{p}}_{p,1})\cap\dot{\mathcal{B}}^{\frac{d}{p}}_{p,1})).$$
(For more details, one can refer to Chapter 3 in \cite{Bahouri2011})
\subsection{Uniqueness}
Assume that $(v_1,F_1,m_1,P_1)$ and $(v_2,F_2,m_2,P_2)$ are two solutions of (\ref{1.5}) with the same initial data. Assume that $\delta v=v_1-v_2,\quad \delta F=F_1-F_2,\quad \delta m=m_1-m_2$. Then we have
 \begin{align}
\left\{
\begin{array}{ll}
\partial_{t}\delta v+(v_1\cdot\nabla)(\delta v)-(1-\omega)\Delta \delta v+\nabla(P_1-P_2)=div~(\sigma_1-\sigma_2)-\delta v\nabla v_2,  \\[1ex]
\partial_{t}\delta F+(v_1\cdot\nabla)\delta F+\partial_\tau \delta F=-\delta F\cdot \nabla v_2+\delta F\cdot \nabla v_2+F_2\cdot \nabla\delta v+\nabla\delta v,  \\[1ex]
\partial_{t}\delta m+(v_1\cdot\nabla)\delta m+\partial_\tau \delta m-\partial^2_s\delta m+(\nabla v_1:\int^s_0S_1)\partial_s\delta m=-\delta v\nabla m_2+H\partial_sm_2,  \\[1ex]
\delta v|_{t=0}=0, \quad \delta F|_{t=0}=0,\quad \delta m|_{t=0}=0, \\[1ex]
\delta F|_{\tau=0}=0, \quad \delta m|_{\tau=0}=0, \quad \delta_m|_{s=\pm\frac{1}{2}}=0,
\end{array}
\right.
\end{align}
where
\[H=\nabla v_1:\int^s_0S_1-\nabla v_2:\int^s_0S_2, \quad \sigma_i=\int^{\frac{1}{2}}_{-\frac{1}{2}}S_ids, \quad S_i=\displaystyle\int^{+\infty}_0m_i\widetilde{\mathscr{S}}(F_i)d\tau,~~ i=1,2.\]

By a similar calculation as in Section 5.3, we have
\begin{align}
\sup_{t\in[0,T]}\|\delta F\|_{L^\infty_\tau(\dot{B}^{-1+\frac{d}{p}}_{p,1})}\leq C_T T^{\frac{1}{2}} \|\delta v\|_{\widetilde{L}^2_T(\dot{B}^{\frac{d}{p}}_{p,1})},
\end{align}
\begin{multline}
\sup_{t\in[0,T]}\|(\lambda+\tau)^2\delta m\|^2_{L^\infty_\tau(\dot{\mathcal{B}}^{-1+\frac{d}{p}}_{p,1})}\leq C_T[\| v^2\|_{\widetilde{L}^2_T(\dot{B}^{\frac{d}{p}}_{p,1})}(\|\delta F\|_{L^\infty_\tau(\dot{B}^{-1+\frac{d}{p}}_{p,1})}\\
+\|(\lambda+\tau)^2\delta m\|^2_{L^\infty_\tau(\dot{\mathcal{B}}^{-1+\frac{d}{p}}_{p,1})}))+\|\delta v\|_{\widetilde{L}^2_T(\dot{B}^{\frac{d}{p}}_{p,1})}],
\end{multline}
\begin{multline}
\|\delta v\|_{\widetilde{L}^\infty_T(\dot{B}^{-1+\frac{d}{p}}_{p,1})\cap\widetilde{L}^2_T(\dot{B}^{\frac{d}{p}}_{p,1})}\leq C_T [T^{\frac{1}{2}}(\|\delta F\|_{L^\infty_\tau(\dot{B}^{-1+\frac{d}{p}}_{p,1})}\\
+\|(\lambda+\tau)^2\delta m\|^2_{L^\infty_\tau(\dot{\mathcal{B}}^{-1+\frac{d}{p}}_{p,1})}))+\| v^2\|_{\widetilde{L}^2_T(\dot{B}^{\frac{d}{p}}_{p,1})}\|\delta v\|_{\widetilde{L}^\infty_T(\dot{B}^{-1+\frac{d}{p}}_{p,1})}.
\end{multline}
If $T_0>0$ is small enough, from the above inequalities we can deduce that $\delta F=0, \delta v=0, \delta m=0$ in $[0,T_0]$ by the argument as in Section 5.3.
Splitting the interval $[0,T]$ into $[0,T_1]$, $[T_1,T_2]$,\,...,\,$[T_k,T_{k+1}]$, where each subinterval satisfies $T_{k+1}-T_k<T_0$. Then we see that $(v_1,F_1, m_1)=(v_2,F_2,m_2)$ in the interval $[0,T_1]$. Since $v_1(T_1)=v_2(T_1),~~F_1(T_1)=F_1(T_1),~~m_1(T_1)=m_2(T_2)$, if follows that $(v_1,F_1, m_1)=(v_2,F_2,m_2)$ in the interval $[T_1,T_2]$. Repeating the argument, we prove the uniqueness.

\section{Global existence}
~~~~~   In this section, we prove that the solution is global in time if the initial velocity is small and $\lambda$ is large.

   \textbf{Proof of Theorem 2.2:}  Assume that the lifespan $T^{*}$ is finite.
   If $T_0>0$ is small enough, by virtue of Lemma \ref{L4.6} with $s=\frac{d}{p}, \rho=\rho_1=2$ and $s=-1+\frac{d}{p}, \rho=\rho_1=1$ with respective, we have
   \begin{multline*}
   \|v(t)\|_{\widetilde{L}^\infty_{T_0}(\dot{B}^{\frac{d}{p}}_{p,1})\cap \widetilde{L}^2_{T_0}(\dot{B}^{\frac{d}{p}+1}_{p,1})}\leq C(\|v_0\|_{\dot{B}^{\frac{d}{p}}_{p,1}}+\|\tau\|_{\widetilde{L}^2_{T_0}(\dot{B}^{\frac{d}{p}}_{p,1})})\\
   \leq C(\|v_0\|_{\dot{B}^{\frac{d}{p}}_{p,1}}+\frac{C}{\lambda}\sup_{t\in[0,T_0]}\|F\|_{L^\infty_\tau(\dot{B}^{\frac{d}{p}}_{p,1})}\|(\lambda+\tau)^2 m\|^2_{L^\infty_\tau(\dot{\mathcal{B}}^{\frac{d}{p}}_{p,1})})\leq M\ep,
   \end{multline*}
   \begin{multline*}
   \|v(t)\|_{\widetilde{L}^\infty_{T_0}(\dot{B}^{-1+\frac{d}{p}}_{p,1})\cap L^1_{T_0}(\dot{B}^{\frac{d}{p}+1}_{p,1})}\leq C(\|v_0\|_{\dot{B}^{-1+\frac{d}{p}}_{p,1}}+\|\tau\|_{\widetilde{L}^1_{T_0}(\dot{B}^{\frac{d}{p}}_{p,1})})\\
   \leq C(\|v_0\|_{\dot{B}^{-1+\frac{d}{p}}_{p,1}}+\frac{C}{\lambda}\sup_{t\in[0,T_0]}\|F\|_{L^\infty_\tau(\dot{B}^{\frac{d}{p}}_{p,1})}\|(\lambda+\tau)^2 m\|^2_{L^\infty_\tau(\dot{\mathcal{B}}^{\frac{d}{p}}_{p,1})})\leq M\ep,
   \end{multline*}
where $M>1$ is a uniform constant. Thus we can define $\overline{T}$  as follows:
   \begin{align}\label{6.1}
   \overline{T}=\sup\bigg\{T:\|v(t)\|_{\widetilde{L}^\infty_{T}(\dot{B}^{\frac{d}{p}}_{p,1})\cap \widetilde{L}^2_{T}(\dot{B}^{\frac{d}{p}+1}_{p,1})}+\|v(t)\|_{\widetilde{L}^\infty_{T}(\dot{B}^{-1+\frac{d}{p}}_{p,1})\cap L^1_{T}(\dot{B}^{\frac{d}{p}+1}_{p,1})} \leq 2M\ep \bigg\}.
   \end{align}
   We need to prove that $\overline{T}=T^{*}$. If $\overline{T}<T^{*}$, by virtue of Proposition \ref{P4.1}, we obtain that
   \begin{align}\label{6.2}
   \sup_{t\in[0,\overline{T}]}\|F\|_{L^\infty_\tau(\dot{B}^{\frac{d}{p}}_{p,1})}\leq CMe^{M\ep}(\|F_0\|_{L^\infty_\tau(\dot{B}^{\frac{d}{p}}_{p,1})}+\ep)\leq C.
   \end{align}

From the equations, we have
\begin{align}\label{6.3}
\partial_t\widetilde{m}+(v\cdot\nabla)\widetilde{m}+\partial_\tau \widetilde{m}-\frac{2}{\lambda+\tau}\widetilde{m}+(\nabla v:\int^s_0S)\partial_s\widetilde{m}-\partial^2_s\widetilde{m}=0,
\end{align}
where $\widetilde{m}=(\lambda+\tau)^2m$.  Since $\frac{1}{\lambda}\leq \ep^2$, it follows from Lemma \ref{L7} that
$$\|\widetilde{m}\|_{L^\infty_{\tau,x,s}}\leq \|\widetilde{m}_0\|_{L^\infty_{\tau,x,s}},$$
for any $t<\overline{T}$. From (\ref{4.11}), we see that
\begin{align}\label{6.4}
&\partial_t\|\dot{\Delta}_j\widetilde{m}\|^p_{L^p_{x,s}}+\partial_\tau\|\dot{\Delta}_j\widetilde{m}\|^p_{L^p_{x,s}}+\int_{\mathbb{R}^d}\int^{\frac{1}{2}}_{-\frac{1}{2}}|\partial_s(\dot{\Delta}_j\widetilde{m})|^2|\dot{\Delta}_j\widetilde{m}|^{p-2}dsdx
-\frac{1}{\lambda}\|\dot{\Delta}_j\widetilde{m}\|^p_{L^p_{x,s}}\\
\nonumber&\lm(\|R_j\|_{L^p_{x,s}}+\|\dot{\Delta}_j[(\nabla v:S)\widetilde{m}]\|_{L^p_{x,s}})\|\dot{\Delta}_j\widetilde{m}\|^{p-1}_{L^p_{x,s}}+\|\dot{\Delta}_j[(\nabla v:\int^s_0S)\widetilde{m}]\|^2_{L^p_{x,s}}\|\dot{\Delta}_j\widetilde{m}\|^{p-2}_{L^p_{x,s}},
\end{align}
where $R_j=[v\cdot \nabla, \dot{\Delta_j}]\widetilde{m}$.
 Note that $m|_{s=\pm\frac{1}{2}}=0$. Using the Poincar\'{e} inequality and the fact that $\frac{1}{\lambda}$ is small enough, and we deduce that
\begin{align}\label{6.5}
&\partial_t\|\dot{\Delta}_j\widetilde{m}\|^p_{L^p_{x,s}}+\partial_\tau\|\dot{\Delta}_j\widetilde{m}\|^p_{L^p_{x,s}}+C\|\dot{\Delta}_j\widetilde{m}\|^p_{L^p_{x,s}}\\
\nonumber&\lm(\|R_j\|_{L^p_{x,s}}+\|\dot{\Delta}_j[(\nabla v:S)\widetilde{m}]\|_{L^p_{x,s}})\|\dot{\Delta}_j\widetilde{m}\|^{p-1}_{L^p_{x,s}}+\|\dot{\Delta}_j[(\nabla v:\int^s_0S)\widetilde{m}]\|^2_{L^p_{x,s}}\|\dot{\Delta}_j\widetilde{m}\|^{p-2}_{L^p_{x,s}},
\end{align}
which leads to
\begin{align}\label{6.6}
&\partial_t\|\dot{\Delta}_j\widetilde{m}\|^2_{L^p_{x,s}}+\partial_\tau\|\dot{\Delta}_j\widetilde{m}\|^2_{L^p_{x,s}}+C\|\dot{\Delta}_j\widetilde{m}\|^2_{L^p_{x,s}}\\
\nonumber&\lm(\|R_j\|_{L^p_{x,s}}+\|\dot{\Delta}_j[(\nabla v:S)\widetilde{m}]\|_{L^p_{x,s}})\|\dot{\Delta}_j\widetilde{m}\|_{L^p_{x,s}}+\|\dot{\Delta}_j[(\nabla v:\int^s_0S)\widetilde{m}]\|^2_{L^p_{x,s}}.
\end{align}

Thanks to the Cauchy-Schwarz inequality, we get
\begin{multline}\label{6.7}
\partial_t\|\dot{\Delta}_j\widetilde{m}\|^2_{L^p_{x,s}}+\partial_\tau\|\dot{\Delta}_j\widetilde{m}\|^2_{L^p_{x,s}}+C\|\dot{\Delta}_j\widetilde{m}\|^2_{L^p_{x,s}}\\
\lm \|R_j\|^2_{L^p_{x,s}}+\|\dot{\Delta}_j[(\nabla v:S)\widetilde{m}]\|^2_{L^p_{x,s}}+\|\dot{\Delta}_j[(\nabla v:\int^s_0S)\widetilde{m}]\|^2_{L^p_{x,s}}.
\end{multline}

By the definition of Besov spaces and Lemma \ref{L2}, we have
\[\|R_j\|_{L^p_{x,s}}\lm 2^{-j\frac{d}{p}}c_j\|u\|_{\dot{B}^{1+\frac{d}{p}}_{p,1}}\|\widetilde{m}\|_{\dot{\mathcal{B}}^{\frac{d}{p}}_{p,1}}, \]
\[\|\dot{\Delta}_j[(\nabla u:S)\widetilde{m}]\|_{L^p_{x,s}}\lm 2^{-j\frac{d}{p}}c_j\|(\nabla u:S)\widetilde{m}\|_{\dot{\mathcal{B}}^{\frac{d}{p}}_{p,1}}\lm 2^{-j\frac{d}{p}}c_j\|u\|_{\dot{B}^{1+\frac{d}{p}}_{p,1}}(\|S\|_{L^\infty_{x,s}}\|\widetilde{m}\|_{\dot{\mathcal{B}}^{\frac{d}{p}}_{p,1}}+\|\widetilde{m}\|_{L^\infty_{x,s}}\|S\|_{\dot{\mathcal{B}}^{\frac{d}{p}}_{p,1}}),\]
where $c_j\in l^1$.

By virtue of Corollary \ref{C2} and the inequality (\ref{6.2}), we have
$$\|S\|_{\dot{\mathcal{B}}^{\frac{d}{p}}_{p,1}}\lm \frac{1}{\lambda}\|\widetilde{m}\|_{L^\infty_\tau(\dot{\mathcal{B}}^{\frac{d}{p}}_{p,1})},\quad  \|S\|_{L^\infty_{x,s}}\lm \frac{1}{\lambda},$$
which implies that
\begin{align}\label{6.8}
\|\dot{\Delta}_j[(\nabla u:S)\widetilde{m}]\|_{L^p_{x,s}}\lm \frac{1}{\lambda}2^{-j\frac{d}{p}}c_j\|u\|_{\dot{B}^{1+\frac{d}{p}}_{p,1}}\|\widetilde{m}\|_{\dot{\mathcal{B}}^{\frac{d}{p}}_{p,1}}.
\end{align}
By the same token, we have
\begin{align}\label{6.9}
\|\dot{\Delta}_j[(\nabla u:\int^s_0S)\widetilde{m}]\|_{L^p_{x,s}}\lm \frac{1}{\lambda}2^{-j\frac{d}{p}}c_j\|u\|_{\dot{B}^{1+\frac{d}{p}}_{p,1}}\|\widetilde{m}\|_{\dot{\mathcal{B}}^{\frac{d}{p}}_{p,1}}.
\end{align}
Plugging (\ref{6.8}) and (\ref{6.9}) into (\ref{6.7}), yields that
\begin{align*}
\partial_t\|\dot{\Delta}_j\widetilde{m}\|^2_{L^p_{x,s}}+\partial_\tau\|\dot{\Delta}_j\widetilde{m}\|^2_{L^p_{x,s}}+C\|\dot{\Delta}_j\widetilde{m}\|^2_{L^p_{x,s}}\lm 2^{-2j\frac{d}{p}}c^2_j\|u\|^2_{\dot{B}^{1+\frac{d}{p}}_{p,1}}\|\widetilde{m}\|^2_{\dot{\mathcal{B}}^{\frac{d}{p}}_{p,1}}.
\end{align*}
Let $f(t,\tau)=\|\dot{\Delta}_j\widetilde{m}\|_{L^p_{x,s}}$ and $h(t,\tau)=2^{-j\frac{d}{p}}c_j\|u\|_{\dot{B}^{1+\frac{d}{p}}_{p,1}}\|\widetilde{m}\|_{\dot{\mathcal{B}}^{\frac{d}{p}}_{p,1}}$.
Changing the variable $z_1=\frac{t+\tau}{2}$ and $z_2=\frac{\tau-t}{2}$, we obtain that
\begin{align*}
\partial_{z_1}\overline{f}^2(z_1,z_2)+C\overline{f}^2(z_1,z_2)\lm \overline{h}^2(z_1,z_2),
\end{align*}
where $\overline{f}(z_1,z_2)=f(t,\tau)$ and $\overline{h}(z_1,z_2)=h(t,\tau)$. From the above inequality we deuce that
\begin{align*}
\partial_{z_1}[e^{Cz_1}\overline{f}^2(z_1,z_2)]\lm e^{Cz_1}\overline{h}^2(z_1,z_2),
\end{align*}
Integrating over $[z_2,z_1]$ with respect to $z_1$, we deduce that
\begin{align*}
e^{Cz_1}\overline{f}^2(z_1,z_2)\lm e^{Cz_2}\overline{f}^2(z_2,z_2)+\int^{z_1}_{z_2}e^{Cz}\overline{h}^2(z,z_2)dz.
\end{align*}
Since $z_2-z_1=-t$, $\overline{f}(z_1,z_2)=f(t,\tau)$ and $\overline{h}(z_1,z_2)=h(t,\tau)$, it follows from the above inequality that
\begin{align*}
f^2(t,\tau)dt'\lm e^{-Ct}f^2(0,\tau)+\int^t_0e^{-C(t-t')}h^2(t',\tau)dt',
\end{align*}
that is
\begin{align*}
\|\dot{\Delta}_j\widetilde{m}\|^2_{L^p_{x,s}}\lm e^{-Ct}\|\dot{\Delta}_j\widetilde{m}_0\|^2_{L^p_{x,s}}+\int^t_0e^{-C(t-t')}2^{-2j\frac{d}{p}}c^2_j\|u\|^2_{\dot{B}^{1+\frac{d}{p}}_{p,1}}\|\widetilde{m}\|^2_{\dot{\mathcal{B}}^{\frac{d}{p}}_{p,1}}dt'
\end{align*}
Multiplying by $2^{2j\frac{d}{p}}$ and taking $l^1$-norm on the both sides of the above inequality, we obtain
\begin{align*}
\sup_{t\in[0,\overline{T}]}e^{Ct}\|\widetilde{m}\|_{L^\infty_\tau(\dot{\mathcal{B}}^{\frac{d}{p}}_{p,1})}\lm \|\widetilde{m}_0\|_{\dot{\mathcal{B}}^{\frac{d}{p}}_{p,1}}+\|u\|_{\widetilde{L}^2_{\overline{T}}(\dot{B}^{1+\frac{d}{p}}_{p,1})}\sup_{t\in[0,\overline{T}]}e^{Ct}\|\widetilde{m}\|_{L^\infty_\tau(\dot{\mathcal{B}}^{\frac{d}{p}}_{p,1})}
\end{align*}
Using the fact that $\|u\|_{\widetilde{L}^2_{\overline{T}}(\dot{B}^{1+\frac{d}{p}}_{p,1})}\leq M\ep$ and choosing $\ep$ small enough, we get
\begin{align}\label{6.10}
\|\widetilde{m}\|_{L^\infty_\tau(\dot{\mathcal{B}}^{\frac{d}{p}}_{p,1})}\lm e^{-Ct}\|\widetilde{m}_0\|_{\dot{\mathcal{B}}^{\frac{d}{p}}_{p,1}},
\end{align}
which implies that
\begin{align}\label{6.11}
\|\widetilde{m}\|_{L^{\rho}_T(L^\infty_\tau(\dot{\mathcal{B}}^{\frac{d}{p}}_{p,1}))}\lm \|\widetilde{m}_0\|_{\dot{\mathcal{B}}^{\frac{d}{p}}_{p,1}},
\end{align}
for any $\rho\in[1,+\infty]$.

For the Navier-Stokes equations, we can write that
\begin{align}\label{6.12}
v=e^{t(1-\omega)\Delta}v_0+\int^t_0 e^{(t-t')(1-\omega)\Delta}(div~\sigma-v\cdot \nabla u-\nabla P) dt'.
\end{align}
 Applying $\dot\Delta_j$ to (\ref{6.12}), we have
\begin{align*}
\dot\Delta_j v=e^{t(1-\omega)\Delta}\dot\Delta_j v_0+\int^t_0 e^{(t-t')(1-\omega)\Delta}(div~\dot\Delta_j\sigma-\dot\Delta_j(v\cdot\nabla v)-\nabla \dot\Delta_j P) dt'.
\end{align*}
By Lemma \ref{L1}, we deduce that
\begin{align*}
\|\dot\Delta_j v\|_{L^p}\leq e^{-tc(1-\omega)2^{2j}}\|\dot\Delta_j v_0\|_{L^p}+C\int^t_0 e^{-(t-t')c(1-\omega)2^{2j}}(\|div\dot\Delta_j\sigma\|_{L^p}+\|\dot\Delta_j(v\cdot \nabla v)\|_{L^p}+\|\nabla \dot\Delta_j P\|_{L^p})dt'.
\end{align*}
Since $-\Delta P=div~(v\cdot\nabla v-div~\sigma)$, it follows that $\|\nabla \dot\Delta_j P\|_{L^p}\lm \|div\dot\Delta_j\sigma\|_{L^p}+\|\dot\Delta_j(v\cdot\nabla v)\|_{L^p}$. Thus
\begin{align}\label{6.13}
\|\dot\Delta_j v\|_{L^p}\leq e^{-tc(1-\omega)2^{2j}}\|\dot\Delta_j v_0\|_{L^p}+C\int^t_0 e^{-(t-t')c(1-\omega)2^{2j}}(\|div\dot\Delta_j\sigma\|_{L^p}+\|\dot\Delta_j(v\cdot \nabla v)\|_{L^p})dt'.
\end{align}
Using Young's inequality, we get
\begin{align*}
\|\dot\Delta_j v\|_{L^p} &\leq \|\dot\Delta_j v_0\|_{L^p}+C(\int^t_0 2^{-2j}(\|div\dot\Delta_j\sigma\|^2_{L^p}+\|\dot\Delta_j(v\cdot \nabla v)\|^2_{L^p})dt')^{\frac{1}{2}} \\
&\leq \|\dot\Delta_j v_0\|_{L^p}+ C(\int^t_0 (\|\dot\Delta_j\sigma\|^2_{L^p}+2^{-2j}\|\dot\Delta_j(v\cdot \nabla v)\|^2_{L^p})dt')^{\frac{1}{2}}.
\end{align*}
Multiplying both sides of the above inequality by $2^{j\frac{d}{p}}$ and taking the $l^1$-norm, we infer that
\begin{align}\label{6.14}
\sup_{t\in[0,\overline{T}]}\|v\|_{\dot B^{\frac{d}{p}}_{p,1}} \leq \|v_0\|_{\dot B^{\frac{d}{p}}_{p,1}}+C[\|\sigma\|_{\widetilde{L}^2_{\overline{T}}(\dot B^{\frac{d}{p}}_{p,1})}+\|v\cdot \nabla v\|_{\widetilde{L}^2_{\overline{T}}(\dot B^{\frac{d}{p}-1}_{p,1})}].
\end{align}

Multiplying both sides of (\ref{6.13}) by $2^{j(\frac{d}{p}+1)}$ and taking the $L^2([0,\overline{T}])$-norm, by a similar calculation, we deduce that
\begin{align}\label{6.15}
(1-\omega)\|v\|_{L^2_{\overline{T}}(\dot B^{\frac{d}{p}+1}_{p,1})}\leq \|v_0\|_{\dot B^{\frac{d}{p}}_{p,r}}+C[\|\sigma\|_{\widetilde{L}^2_{\overline{T}}(\dot B^{\frac{d}{p}}_{p,1})}+\|v\cdot \nabla v\|_{\widetilde{L}^2_{\overline{T}}(\dot B^{\frac{d}{p}-1}_{p,1})}].
\end{align}
By virtue of Corollary \ref{C2} and the estimates (\ref{6.2}), (\ref{6.11}), we obtain that
\begin{align}\label{6.16}
\|\sigma\|_{\widetilde{L}^2_{\overline{T}}(\dot B^{\frac{d}{p}}_{p,1})}\lm\|S\|_{\widetilde{L}^2_{\overline{T}}(\dot{\mathcal{B}}^{\frac{d}{p}}_{p,1})} \lm\frac{1}{\lambda}\|\widetilde{m}\|_{\widetilde{L}^2_{\overline{T}}(L^\infty_\tau(\dot{\mathcal{B}}^{\frac{d}{p}}_{p,1}))} \sup_{t\in[0,\overline{T}]}\|F\|_{L^\infty_\tau(\dot{B}^{\frac{d}{p}}_{p,1})}\lm \ep^2.
\end{align}
Taking advantage of Corollary \ref{C1}, we see that
\begin{align}\label{6.17}
\|v\cdot \nabla v\|_{\widetilde{L}^2_{\overline{T}}(\dot B^{-1+\frac{d}{p}}_{p,1})}\lm \|v\|_{L^\infty_{\overline{T}}(\dot B^{-1+\frac{d}{p}}_{p,1})}\|\nabla v\|_{\widetilde{L}^2_{\overline{T}}(\dot B^{\frac{d}{p}}_{p,1})}\lm \|v\|_{L^\infty_{\overline{T}}(\dot B^{-1+\frac{d}{p}}_{p,1})}\| v\|_{\widetilde{L}^2_{\overline{T}}(\dot B^{\frac{d}{p}+1}_{p,1})}\lm \ep^2.
\end{align}
Plugging (\ref{6.16}), (\ref{6.17}) into (\ref{6.15}), (\ref{6.16}) yields that
\begin{align}
\|v(t)\|_{\widetilde{L}^\infty_{\overline{T}}(\dot{B}^{\frac{d}{p}}_{p,1})\cap \widetilde{L}^2_{T}(\dot{B}^{\frac{d}{p}+1}_{p,1})} \leq \ep+C\ep^2.
\end{align}
If $\ep$ is small enough, such that $\ep+C\ep^2<M\ep$, then we see that \[\|v(t)\|_{\widetilde{L}^\infty_{\overline{T}}(\dot{B}^{\frac{d}{p}}_{p,1})\cap \widetilde{L}^2_{T}(\dot{B}^{\frac{d}{p}+1}_{p,1})}<M\ep\].

On the other hand, by virtue of (\ref{6.13}), we have
\begin{align*}
\|\dot\Delta_j v\|_{L^p} &\leq \|\dot\Delta_j v_0\|_{L^p}+(\int^t_0 (\|div\dot\Delta_j\sigma\|_{L^p}+\|\dot\Delta_j(v\cdot \nabla v)\|_{L^p})dt') \\
&\leq 2^j\|\dot\Delta_j v_0\|_{L^p}+ C(\int^t_0 (\|\dot\Delta_j\sigma\|_{L^p}+\|\dot\Delta_j(v\cdot \nabla v)\|_{L^p})dt').
\end{align*}
Multiplying both sides of the above inequality by $2^{j(-1+\frac{d}{p})}$ and taking the $l^1$-norm, we infer that
\begin{align}\label{6.19}
\sup_{t\in[0,\overline{T}]}\|v\|_{\dot B^{-1+\frac{d}{p}}_{p,1}} \leq \|v_0\|_{\dot B^{-1+\frac{d}{p}}_{p,1}}+C[\|\sigma\|_{L^1_{\overline{T}}(\dot B^{\frac{d}{p}}_{p,1})}+\|v\cdot \nabla v\|_{L^1_{\overline{T}}(\dot B^{\frac{d}{p}-1}_{p,1})}].
\end{align}
Multiplying both sides of (\ref{6.13}) by $2^{j(\frac{d}{p}+1)}$ and taking the $L^1([0,\overline{T}])$-norm, we deduce that
\begin{align}\label{6.20}
(1-\omega)\|v\|_{L^1_{\overline{T}}(\dot B^{\frac{d}{p}+1}_{p,1})}\leq \|v_0\|_{\dot B^{-1+\frac{d}{p}}_{p,r}}+C[\|\sigma\|_{L^1_{\overline{T}}(\dot B^{\frac{d}{p}}_{p,1})}+\|v\cdot \nabla v\|_{L^1_{\overline{T}}(\dot B^{\frac{d}{p}-1}_{p,1})}].
\end{align}
By virtue of Corollary \ref{C2} and the estimates (\ref{6.2}), (\ref{6.11}), we obtain that
\begin{align}\label{6.21}
\|\sigma\|_{L^1_{\overline{T}}(\dot B^{\frac{d}{p}}_{p,1})}\lm\|S\|_{L^1_{\overline{T}}(\dot{\mathcal{B}}^{\frac{d}{p}}_{p,1})} \lm\frac{1}{\lambda}\|\widetilde{m}\|_{L^1_{\overline{T}}(L^\infty_\tau(\dot{\mathcal{B}}^{\frac{d}{p}}_{p,1}))} \sup_{t\in[0,\overline{T}]}\|F\|_{L^\infty_\tau(\dot{B}^{\frac{d}{p}}_{p,1})}\lm \ep^2.
\end{align}
Taking advantage of Corollary \ref{C1}, we see that
\begin{align}\label{6.22}
\|v\cdot \nabla v\|_{\widetilde{L}^1_{\overline{T}}(\dot B^{-1+\frac{d}{p}}_{p,1})}\lm \|v\|_{L^\infty_{\overline{T}}(\dot B^{-1+\frac{d}{p}}_{p,1})}\|\nabla v\|_{\widetilde{L}^1_{\overline{T}}(\dot B^{\frac{d}{p}}_{p,1})}\lm \|v\|_{L^\infty_{\overline{T}}(\dot B^{-1+\frac{d}{p}}_{p,1})}\| v\|_{\widetilde{L}^1_{\overline{T}}(\dot B^{\frac{d}{p}+1}_{p,1})}\lm \ep^2.
\end{align}
Plugging (\ref{6.22}), (\ref{6.21}) into (\ref{6.19}), (\ref{6.20}) yields that
\begin{align}
\|v(t)\|_{\widetilde{L}^\infty_{\overline{T}}(\dot{B}^{-1+\frac{d}{p}}_{p,1})\cap L^1_{T}(\dot{B}^{\frac{d}{p}+1}_{p,1})} \leq \ep+C\ep^2.
\end{align}
If $\ep$ is small enough, such that $\ep+C\ep^2<M\ep$, then we see that \[\|v(t)\|_{\widetilde{L}^\infty_{\overline{T}}(\dot{B}^{-1+\frac{d}{p}}_{p,1})\cap L^1_{T}(\dot{B}^{\frac{d}{p}+1}_{p,1})}<M\ep.\]

 This is a contradiction because of the definition of $\overline{T}$. Thus, we get $\overline{T}=T^*$. By virtue of the estimates (\ref{6.2}), (\ref{6.10}), we verify that
\[\limsup_{t\rightarrow T^{*}}\|F\|_{L^\infty_\tau(\dot{B}^{\frac{d}{p}}_{p,1})}<+\infty, \limsup_{t\rightarrow T^{*}} \|\widetilde{m}\|_{L^\infty_\tau(\dot{\mathcal{B}}^{\frac{d}{p}}_{p,1})}<+\infty.\]

   This implies the solution can be extended beyond $[0,T^{*})$, which contradicts the lifespan $T^{*}$. Thus the solution is global. This completes the proof.\\

\begin{rema}
From the estimate \ref{6.10}, we see that the Besov norm of $m$ is exponent decay in time under the small condition of initial data and $\lambda$ sufficient large.
\end{rema}

\smallskip
\noindent\textbf{Acknowledgments} This work was
partially supported by the National Natural Science Foundation of China (No.11671407 and No.11701586), the Macao Science and Technology Development Fund (No. 098/2013/A3), and Guangdong Province of China Special Support Program (No. 8-2015),
and the key project of the Natural Science Foundation of Guangdong province (No. 2016A030311004). 

\phantomsection
\addcontentsline{toc}{section}{\refname}
\bibliographystyle{abbrv} 
\bibliography{Feneref}

\end{document}